\documentclass[journal,10pt,onecolumn,draftclsnofoot,]{IEEEtran}

\usepackage{amsmath, amssymb, amsthm, mathtools}    % Math packages
\usepackage{adjustbox}      % To adjust size of pictures
\usepackage{standalone}     % To add standalone tikz pictures
\usepackage{xcolor}       % To make color marks in the text
\usepackage{tikz, pgfplots}   % To draw vector graphics and include matlab2tikz figures
\usepackage{float}        % To set figures position (look: http://tex.stackexchange.com/questions/8625/force-figure-placement-in-text)
\usepackage{caption}        % To centering figure captions
\usepackage{subcaption}
\usepackage[shortlabels]{enumitem}  % To customize description environment
\usepackage{verbatim}     % To add comments
\usepackage{algorithm}      % To write algorithms in pseudo code
\usepackage[noend]{algpseudocode} %No end to remove endif, endfor etc.
\usepackage{caption}      % to add title to table
\usepackage{array}        % to specify width of table cell
\usepackage{bm}
\usepackage{todonotes}
\usepackage{xcolor}
\usepackage{subfiles}
\usepackage{pifont}       % for checkmarks and and xmark
\newcommand{\RNum}[1]{\uppercase\expandafter{\romannumeral #1\relax}}
\usetikzlibrary{arrows.meta}    % Arrows library
\usetikzlibrary{shapes.geometric,calc,patterns,angles,quotes} % Geometric library
\usepgfplotslibrary{fillbetween}
\usetikzlibrary{patterns}
\usetikzlibrary{calc}
\pgfplotsset{compat=1.11}

\usepgfplotslibrary{groupplots}

\theoremstyle{plain}

\newtheorem{lemma}{Lemma}

\newtheorem{proposition}{Proposition}
\newtheorem*{proposition*}{Proposition}

\theoremstyle{definition}

\theoremstyle{remark}
\newtheorem*{remark}{Remark}
\newtheorem*{corollary}{Corollary}
\DeclareMathOperator{\argmin}{argmin}
\DeclareMathOperator{\minimise}{minimise}
\DeclareMathOperator{\st}{subject ~ to}
\DeclareMathOperator{\diag}{\mathbf{diag}}

\DeclareMathOperator{\image}{\mathbf{im}}
\DeclareMathOperator{\Real}{Re}
\DeclareMathOperator{\Imag}{Im}

\algdef{SE}[DOWHILE]{Do}{doWhile}{\algorithmicdo}[1]{\algorithmicwhile\ #1}%
\presetkeys{todonotes}{color=yellow!20}{}

\def\BibTeX{{\rm B\kern-.05em{\sc i\kern-.025em b}\kern-.08emT\kern-.1667em\lower.7ex\hbox{E}\kern-.125emX}}
    
\begin{document}

%
% The "title" command has an optional parameter, allowing the author to define a "short title" to be used in page headers.
\title{Low-Voltage Distribution Network Impedance Identification Based on Smart Meter Data}

\author{Sergey~Iakovlev,
        Robin~J.~Evans,~\IEEEmembership{Life Fellow,~IEEE,}
        Iven~Mareels,~\IEEEmembership{Fellow,~IEEE,}
       % Girish Nair,~\IEEEmembership{Member,~IEEE,}% <-this % stops a space
%\thanks{M. Shell was with the Department
%of Electrical and Computer Engineering, Georgia Institute of Technology, Atlanta,
%GA, 30332 USA e-mail: (see http://www.michaelshell.org/contact.html).}% <-this % stops a space
%\thanks{J. Doe and J. Doe are with Anonymous University.}% <-this % stops a space
%\thanks{Manuscript received April 19, 2005; revised August 26, 2015.}}
\thanks{S. Iakovlev and R. J. Evans are with the Department of Electrical
and Electronic Engineering, The University of Melbourne, Melbourne, VIC 3010,
Australia. (e-mail: siakovlev@unimelb.edu.au)}
\thanks{I. Mareels is with IBM Research, Melbourne, VIC 3010,
Australia}}

\maketitle

% The abstract is a short summary of the work to be presented in the article.
\begin{abstract}
Under conditions of high penetration of renewables, the low-voltage (LV) distribution network needs to be carefully managed. In such a scenario, an accurate real-time low-voltage power network model is an important prerequisite, which opens up the possibility for application of many advanced network control and optimisation methods thus providing improved power flow balancing, reduced maintenance costs, and enhanced reliability and security of a grid. 

Smart meters serve as a source of information in LV networks and allow for accurate measurements at almost every node, which makes it advantageous to use data driven methods. In this paper, we formulate a non-linear and non-convex problem, solve it efficiently, and propose a number of fully smart meter data driven methods for line parameters estimation. Our algorithms are fast, recursive in data, scale linearly with the number of nodes, and can be executed in a decentralised manner. The performance of these algorithms is demonstrated for different measurement accuracy scenarios through simulations.
\end{abstract}

\begin{IEEEkeywords}
low-voltage distribution grid, smart meter measurements, impedances identification
\end{IEEEkeywords}

\section{Introduction}
The deployment of smart meters at the level of a single customer in the LV power grid opens up the possibility of advanced control and monitoring functions, identification of faults and detection of topology changes. However the most important and developed techniques in state estimation \cite{huang2012state, della2014electrical}, optimal power flow \cite{lam2012distributed}, active filtering \cite{tarkiainen2004identification} and economic dispatch all require line impedances and the topology of feeder systems to be known \cite{wood2012power, jahangiri2013distributed}. This information is not always available and when it is available it is often inaccurate. Moreover the increasing emergence of plug-and-play parts in the modern LV distribution grid (e.g. hybrid electric vehicles, renewable generators etc.) further exacerbates this situation. Additionally, certain assumptions for transmission lines in the high voltage (HV) network do not hold for the LV case. For example PMU measurements are a mature technology for HV grids, whereas in LV distribution grids PMU measurements are usually not available primarily due to the high cost \cite{de2010synchronized}. However, the analysis of LV distribution grid estimation in the existing literature often assumes the availability of PMU synchronised measurements \cite{yang2010online, cavraro2015data, deka2016estimating}. These factors make it difficult to apply the control and optimisation methods mentioned above. Hence building a LV grid model based on SM data is an important and potentially valuable opportunity. Smart meters serve as a source of data and provide reliable information about the LV grid which can enable realistic modelling \cite{alahakoon2016smart}. 

The most recent research activity is focused on impedance identification methods that allow identification of every single power line impedance (\cite{han2016automated, yang2010online}) rather than the grid equivalent impedance as previously \cite{cobreces2009grid, ciobotaru2011line}. With accurate models these methods open up the possibility for real-time tracking of aging related degradation, faults and electricity theft detection and localisation. Theft is indeed a real issue, e.g. 50\% of electricity in developing countries \cite{tariq2016electricity}, 1-3\% (\$6 billion) of total revenue in US \cite{van2011electricity} and 1\% (1200GWh) of electrical energy in Netherlands annually \cite{sahoo2015electricity}. The most common approach to stealing electricity is a direct connection to the low voltage grid bypassing meter infrastructure completely. Real-time estimation of single power line impedances can help to localise such situations \cite{tariq2016electricity}.  

This paper presents four contributions. First, we propose a recursive approach to low-voltage network modelling based on smart meter data only. Second, we provide a fully decentralised scalable system identification method based on this model and prove its optimality under certain conditions. Third, the method we propose finds a global solution for a class of non-convex optimisation problems in iterative fashion. Fourth, we consider two practical modifications that can improve algorithm performance for industrial applications. Importantly, our algorithms do not require PMU synchronised measurements as is often assumed \cite{de2010synchronized}. 

The paper is organised as follows. Section \RNum{2} introduces a recursive model for the low-voltage power distribution grid, notations, and presents the main assumptions used through the paper. Section \RNum{3} provides a description of the proposed identification approach and introduces the idea of decentralised grid identification. In Section \RNum{4}, the algorithm is tested via MATLAB simulations on IEEE test feeder data and its applicability and performance are illustrated. Section \RNum{5} concludes the paper. The appendix contains proofs of certain theoretical results that are used in the paper. 

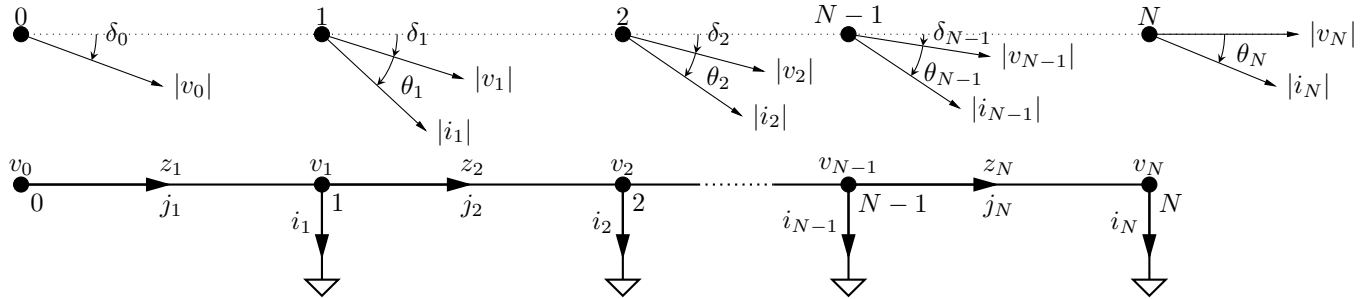
\begin{figure*}[!t]
  %\captionsetup{justification=centering}
  \centering
  \begin{tikzpicture}
    \draw[fill] (0,0) circle [radius=0.1] node[above]{$v_{0}$} node[below right]{$0$};
    \draw[thick] (0,0) -- (9,0) node[pos=2/9, above]{$z_{1}$} node[pos=6/9, above]{$z_{2}$};
    \draw[arrows = {-Stealth[inset=0pt, angle=30:10pt]}, thick] (0,0) -- (2,0)  node[below]{$j_{1}$};
    \draw[fill] (4,0) circle [radius=0.1] node[above]{$v_{1}$} node[below right]{$1$};
    \draw[arrows = {-Stealth[inset=0pt, angle=30:10pt]}, thick] (4,0) -- (6,0)  node[below]{$j_{2}$};
    \draw[arrows = {-Stealth[fill=none, inset=0pt, angle=90:10pt]}, thick] (4,0) -- (4,-1.5);
    \draw[arrows = {-Stealth[inset=0pt, angle=30:10pt]}, thick] (4,0) -- node[left]{$i_{1}$}(4,-1);
    \draw[fill] (8,0) circle [radius=0.1] node[above]{$v_{2}$} node[below right]{$2$};
    \draw[arrows = {-Stealth[fill=none, inset=0pt, angle=90:10pt]}, thick] (8,0) -- (8,-1.5);
    \draw[arrows = {-Stealth[inset=0pt, angle=30:10pt]}, thick] (8,0) -- node[left]{$i_{2}$}(8,-1);
    \draw[dotted, thick] (9,0) -- (10,0);
    \draw[thick] (10,0) -- (11,0);
    \draw[fill] (11,0) circle [radius=0.1] node[above]{$v_{N-1}$} node[below right]{$N-1$};
    \draw[arrows = {-Stealth[fill=none, inset=0pt, angle=90:10pt]}, thick] (11,0) -- (11,-1.5);
    \draw[arrows = {-Stealth[inset=0pt, angle=30:10pt]}, thick] (11,0) -- node[left]{$i_{N-1}$}(11,-1);
    \draw[thick] (11,0) -- node[above]{$z_{N}$}(15,0);
    \draw[fill] (15,0) circle [radius=0.1] node[above]{$v_{N}$} node[below right]{$N$};
    \draw[arrows = {-Stealth[inset=0pt, angle=30:10pt]}, thick] (11,0) -- (13,0)  node[below]{$j_{N}$};
    \draw[arrows = {-Stealth[fill=none, inset=0pt, angle=90:10pt]}, thick] (15,0) -- (15,-1.5);
    \draw[arrows = {-Stealth[inset=0pt, angle=30:10pt]}, thick] (15,0) -- node[left]{$i_{N}$}(15,-1);

    \draw[dotted, thin] (0, 2) -- (17,2);
    \draw[fill] (0,2) circle [radius=0.1] node[above]{$0$};
    \draw[fill] (4,2) circle [radius=0.1] node[above]{$1$};
    \draw[fill] (8,2) circle [radius=0.1] node[above]{$2$};
    \draw[fill] (11,2) circle [radius=0.1] node[above]{$N-1$};
    \draw[fill] (15,2) circle [radius=0.1] node[above]{$N$};

    \coordinate (a1) at (0,2);
    \coordinate (b1) at (2,2);
    \coordinate (c1) at (1.9,1.3);
    \draw[arrows = {-Stealth[inset=0pt, angle=25:5pt]}] (a1) node[right = 1cm]{$\delta_{0}$}--(c1) node[right]{$|v_0|$};
    \pic[draw,<-,>=stealth,angle eccentricity=1.3, angle radius=1cm] {angle=c1--a1--b1};

    \coordinate (a2) at (4,2);
    \coordinate (b2) at (6,2);
    \coordinate (c2) at (5.9,1.4);
    \coordinate (d2) at (5.4,0.7);
    \draw[arrows = {-Stealth[inset=0pt, angle=25:5pt]}] (a2) node[right = 1cm]{$\delta_{1}$}--(c2) node[right]{$|v_1|$};
    \pic[draw,<-,>=stealth,angle eccentricity=1.3, angle radius=1cm] {angle=c2--a2--b2};
    \draw[arrows = {-Stealth[inset=0pt, angle=25:5pt]}] (a2) -- (d2) node[right]{$|i_1|$};
    \pic["$\theta_1$", draw,<-,>=stealth,angle eccentricity=1.4, angle radius=1cm] {angle=d2--a2--c2};

    \coordinate (a3) at (8,2);
    \coordinate (b3) at (10,2);
    \coordinate (c3) at (9.9,1.5);
    \coordinate (d3) at (9.6,0.9);
    \draw[arrows = {-Stealth[inset=0pt, angle=25:5pt]}] (a3) node[right = 1cm]{$\delta_{2}$}--(c3) node[right]{$|v_2|$};
    \pic[draw,<-,>=stealth,angle eccentricity=1.3, angle radius=1cm] {angle=c3--a3--b3};
    \draw[arrows = {-Stealth[inset=0pt, angle=25:5pt]}] (a3) -- (d3) node[right]{$|i_2|$};
    \pic["$\theta_2$", draw,<-,>=stealth,angle eccentricity=1.4, angle radius=1cm] {angle=d3--a3--c3};

    \coordinate (aN-1) at (11,2);
    \coordinate (bN-1) at (13,2);
    \coordinate (cN-1) at (12.9,1.7);
    \coordinate (dN-1) at (12.5,1);
    \draw[arrows = {-Stealth[inset=0pt, angle=25:5pt]}] (aN-1) node[right = 1cm]{$\delta_{N-1}$}--(cN-1) node[right]{$|v_{N-1}|$};
    \pic[draw,<-,>=stealth,angle eccentricity=1.6, angle radius=1cm] {angle=cN-1--aN-1--bN-1};
    \draw[arrows = {-Stealth[inset=0pt, angle=25:5pt]}] (aN-1) -- (dN-1) node[right]{$|i_{N-1}|$};
    \pic["$\theta_{N-1}$", draw,<-,>=stealth,angle eccentricity=1.5, angle radius=1cm] {angle=dN-1--aN-1--cN-1};

    \coordinate (aN) at (15,2);
    \coordinate (bN) at (17,2);
    \coordinate (cN) at (16.7,1.3);
    \draw[arrows = {-Stealth[inset=0pt, angle=25:5pt]}] (aN)--(bN) node[right]{$|v_N|$};
    \draw[arrows = {-Stealth[inset=0pt, angle=25:5pt]}] (aN) -- (cN) node[right]{$|i_N|$};
    \pic["$\theta_N$", draw,<-,>=stealth,angle eccentricity=1.4, angle radius=1cm] {angle=cN--aN--bN};

  \end{tikzpicture}
  \caption{Bottom: model of a single-phase chain low-voltage distribution feeder. 
  Top: global voltage ($v^g_n$) and current ($i_{n}e^{i\delta_n}$) phasors at each node. In this figure the phase reference is chosen with respect to the last node $N$.}
  \label{f1}
\end{figure*}

\section{Model Formulation}
In this section we review a circuit theory approach for low voltage power network analysis by considering a chain feeder first and then generalise the results for a tree network.

\subsection{Notation and preliminaries}
The following conventions are used for description of a low-voltage network:
\begin{itemize}
  \item Node $0$ represents a substation transformer.
  \item Chain feeder nodes and lines are indexed such that $n$-th power line connects nodes $n-1$ and $n$, we denote it as $(n-1, n)$. Therefore all chain feeder parameters have one index.
  \item A radial power network is represented by a tree graph $G(\mathcal{N}, \mathcal{E})$ where each node in $\mathcal{N}$ denotes a bus number and each link $(k, n) \in \mathcal{E}$ denotes a power line between nodes $k$ and $n$. Power line related parameters in the radial network require two indices.
  \item Small Latin letters are used for complex or real scalars, i.e. $a_n = |a_n|e^{i\angle a_n}$ where $|a_n|$ - amplitude of the complex scalar $a_n$ and $\angle a_n$ - its angle. Bold letters are used for vectors ($\bm{a}_n = \big[a_{1, n}, \ldots, a_{M, n}\big]^T = \big[a_{1:M, n}\big]^T$) and bold capital letters for matrices ($\bm{A}_n = \big[a_{1:M, 1:N}\big] = \big[\bm{a}_{1:N}\big]$). 
\end{itemize}
We limit our study to steady state behaviour of a single phase LV grid, when all voltages and currents can be represented as phasors.

\subsection{Chain feeder. Models} 

Consider a single phase model of a chain low-voltage distribution feeder (Fig. \ref{f1}) where $N$ buses (nodes $1, \ldots, N$) are connected in series to a distribution transformer (node $0$). Every bus $n$ has voltage $v_n$ and current $i_n$, that are represented as phasors with angle $\theta_n$ between them. When $\theta_n > 0$, current is lagging the voltage. Power line impedances and corresponding line currents are denoted by $z_n$ and $j_n$ respectively. By $\delta_n$ we denote the phase of a voltage phasor with respect to a global reference, for example in Figure \ref{f1} the phase reference is chosen with respect to the last node $N$.

Throughout this paper we distinguish between global variables, i.e. variables defined with respect to a global phase reference, and local variables, i.e. variables that do not contain global phase angle ($\delta_n$) information:

\begin{itemize}
  \item $v_n \coloneqq |v_n|$ - local voltage at the node $n$;
  \item $v_n^g \coloneqq |v_n|e^{i\delta_n}$ - global voltage at the node $n$;
  \item $i_{n} \coloneqq |i_{n}|e^{i\theta_n}$ - local current (consumed or injected) at the node $n$;
  \item $i^g_{n} \coloneqq |i_{n}|e^{i(\theta_n + \delta_n)} = i_{n}e^{i\delta_n}$ - global current (consumed or injected) at the node $n$;
\end{itemize}

Given only local information $v_n, i_n, \theta_n, z_n$ for all nodes $n$, we can calculate unknown node phases $\delta_n$ and line currents $j_n$ using Propositions \ref{backward_prop} and \ref{forward-prop}. 

\begin{proposition}[Backward model]
  Refer to Fig. \ref{f1}. Let the phase reference be chosen with respect to the last node $N$ (as in Fig. \ref{f1}), i.e. $\delta_N = 0$. Also, let $j_{n} \coloneqq |j_{n}|e^{i\beta_n}$ is a local line current, then
  \begin{equation}
    \begin{split}
    v_{n-1}e^{i\Delta_{n}} - v_{n} = j_{n}z_{n}, \\
    j_{n-1} = i_{n-1} + j_{n}e^{-i\Delta_{n}}
    \end{split}
    \label{l_pi_original}
  \end{equation}
  for all $n$, where $\Delta_{n} = \delta_{n-1} - \delta_{n}$ is the phase increment corresponding to the power line $(n-1, n)$. 
  \label{backward_prop}
\end{proposition}

\begin{proof}
Define a global line current, flowing between nodes $n-1$ and $n$ as $j^{g}_{n} \coloneqq j_{n}e^{i\delta_n} = |j_{n}|e^{i(\beta_n + \delta_n)}$, i.e. its phase is calculated with respect to the (receiving) node $n$, which gives $j^{g}_{N} = i_N$. Using Ohm's law, we can find the ($N-1$)-th node's phase $\delta_{N-1}$: 
  $$v_{N-1}e^{i\delta_{N-1}} - v_{N} = j^{g}_{N}z_{N},$$
and from Kirchoff's current law we calculate current $j^{g}_{N-1} = i^g_{N-1} + j^{g}_{N}$ flowing between nodes $N-1$ and $N$.
We repeat this procedure for nodes $N-2, N-3, \ldots, 0$, find all phases and obtain the following recursive relations for phases and line currents in the feeder:
\begin{equation*}
  \begin{split}
  v_{n-1}e^{i\delta_{n-1}} - v_{n}e^{i\delta_n} = j^{g}_{n}z_{n}, \\
  j^{g}_{n-1} = i^g_{n-1} + j^{g}_{n}.
  \end{split}
\end{equation*}
Divide the first equation by $e^{i\delta_{n}}$, the second by $e^{i\delta_{n-1}}$, and then noting that $j^{g}_n e^{-i\delta_n} = j_n$, we obtain the model (\ref{l_pi_original}) written in terms of phase increments and local variables.
\end{proof}

It is often convenient to choose the phase reference with respect to the distribution transformer which leads to the alternative model. 
\begin{proposition}[Forward model]
  Refer to Fig. \ref{f1}. Let the phase reference be chosen with respect to the first node in the network, i.e. $\delta_0 = 0$. Also, let $j_{n} \coloneqq |j_{n}|e^{i\phi_n}$ is a local line current, then
  \begin{equation}
    \begin{split}
    v_{n-1} - v_{n}e^{-i\Delta_{n}} = j_{n}z_{n}, \\
    j_{n+1} = j_{n}e^{i\Delta_{n}} - i_{n}
    \end{split}
    \label{lf-forward1}
  \end{equation}
  for all $n$, where $\Delta_{n} = \delta_{n-1} - \delta_{n}$ is the phase increment corresponding to the power line $(n-1, n)$. 
  \label{forward-prop}
\end{proposition}

\begin{proof}
  Define a global line current as $j^{g}_{n} \coloneqq j_{n}e^{i\delta_{n-1}} = |j_{n}|$ $e^{i(\phi_n  \\ + \delta_{n-1})}$, i.e. its phase is defined with respect to the (sending) node $n-1$. This definition enables us to start consideration from the first node of the chain feeder, where $\delta_0 = 0$. Following the same procedure as in the proof of Proposition \ref{backward_prop} we obtain:
  \begin{equation*}
    \begin{split}
    v_{n-1}e^{i\delta_{n-1}} - v_{n}e^{i\delta_n} = j^{g}_{n}z_{n}, \\
    j^{g}_{n+1} = j^{g}_{n} - i^g_{n}.
    \end{split}
  \end{equation*}
  Divide the first equation by $\delta_{n-1}$ and the second one by $\delta_{n}$ we get a model (\ref{lf-forward1}) written in terms of phase increments and local variables. 
\end{proof}

\begin{remark} \hfill
\begin{itemize}
  \begin{comment}
  \item Another useful forward model formulation is obtained after excluding the phase term $e^{i\Delta_{n}}$:
        \begin{equation}
        \begin{split}
        |v_{n-1}| - |v_{n}|e^{-i\Delta_{n}} = j_{n}z_{n} \\
        j_{n+1} = j_{n}\frac{|v_{n-1}|}{|v_{n}|} - z_{n}^*\frac{|j_{n}|^2}{|v_{n}|} - i_{n}
        \end{split}
        \label{lf-forward2}
        \end{equation}
  \end{comment}
  \item The forward model has some limitations. For instance, it requires knowledge of the substation transformer current (i.e. $j_{1}$ on fig. \ref{f1}). In the context of impedance identification this means that its magnitude and phase should be known from measurements. However, this is not always the case for LV distribution grids. Another limitation will become apparent when considering the tree network case. 
  \item The forward and backward models are essentially equivalent in that the forward model (\ref{lf-forward1}) can be obtained from (\ref{l_pi_original}) by multiplying the first equation by $e^{-i\Delta_{n}}$ and noting the difference in phases for line current definitions.  
\end{itemize} 
\end{remark}

\subsection{Tree network}{}

The results of the previous subsection can be extended to the case of a tree topology. We start from a simple branching example, show how to resolve it using the system of equations (\ref{l_pi_original}) and then generalise this approach for any tree network.
\begin{proposition}
  The backward model (\ref{l_pi_original}) resolves any network with a tree topology.
\end{proposition}{}
\begin{proof}[Proof (Outline)]
  Consider the network topology case depicted on Fig. \ref{f2} which consists of three buses{} connected in the form of the simplest branching.
  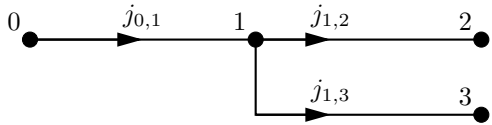
\begin{figure}[H]
  \centering
  \begin{tikzpicture}
  \draw[fill] (0, 0) circle [radius=0.1] node[above left]{$0$};
  \draw[thick] (0, 0) -- (3, 0);
  \draw[arrows = {-Stealth[inset=0pt, angle=30:10pt]}, thick] (0,0) -- (1.5,0)  node[above]{$j_{0,1}$};
  \draw[fill] (3, 0) circle [radius=0.1] node[above left]{$1$};
  %\draw[arrows = {-Stealth[fill=none, inset=0pt, angle=90:10pt]}, thick] (3, 0) -- (3, -2);
  %\draw[arrows = {-Stealth[inset=0pt, angle=30:10pt]}, thick] (3,0) -- (3, -1)  node[right]{$i_{1}$};
  \draw[thick] (3, 0) -- (6, 0);
  \draw[fill] (6, 0) circle [radius=0.1] node[above left]{$2$};
  %\draw[arrows = {-Stealth[fill=none, inset=0pt, angle=90:10pt]}, thick] (6,0) -- (6,-2);
  %\draw[arrows = {-Stealth[inset=0pt, angle=30:10pt]}, thick] (6,0) -- (6, -1)  node[right]{$i_{2}$};
  \draw[arrows = {-Stealth[inset=0pt, angle=30:10pt]}, thick] (3,0) -- (4, 0)  node[above]{$j_{1,2}$};
  \draw[thick] (3, 0) -- (3, -1) -- (6, -1);
  \draw[fill] (6, -1) circle [radius=0.1] node[above left]{$3$};
  %\draw[thick] (3, 2) -- (2, 2);
  %\draw[arrows = {-Stealth[inset=0pt, angle=30:10pt]}, ultra thick, dotted] (6.5,1) -- (5,1)  node[pos=1/2, above]{\textbf{BCI}$_1$};
  %\draw[arrows = {-Stealth[inset=0pt, angle=30:10pt]}, ultra thick, dotted] (4,2.5) -- (4,1)  node[pos=1/2, right]{\textbf{BCI}$_2$};
  %\draw[arrows = {-Stealth[inset=0pt, angle=30:10pt]}, thick] (3, 2) -- (2, 2)  node[above]{$i_{3}$};
  %\draw[arrows = {-Stealth[fill=none, inset=0pt, angle=90:10pt]}, thick] (2, 2) -- (2, 1);
  \draw[arrows = {-Stealth[inset=0pt, angle=30:10pt]}, thick] (3, -1) -- (4, -1) node[above]{$j_{1,3}$};
  %\draw[dotted, thick] (6, 0) -- (6.5, 0);
  %\draw[dotted, thick] (3, 2) -- (3, 2.5);
  \end{tikzpicture}
  \caption{The simplest branching case}
  \label{f2}
  \end{figure}
  Apply backward model (\ref{l_pi_original}) to branches $(1, 3)$ and $(1, 2)$ simultaneously and find $\Delta_{1,3}$, $\Delta_{1,2}$, $j_{0, 1}$:
  \begin{equation}
    \begin{split}
      &{} v_{1}e^{i\Delta_{1,3}} - v_{3} = j_{1,3}z_{1,3}, \\
      & v_{1}e^{i\Delta_{1,2}} - v_{2} = j_{1,2}z_{1,2},
    \end{split}
  \end{equation}
  where $j_{1,3} = i_3$ and $j_{1,2} = i_2$. Then 
  $$j_{0, 1} = j_{0, 1}^ge^{-i\delta_{1}} = j_{1,3}e^{-i\Delta_{1,3}} + j_{1,2}e^{-i\Delta_{1,2}} + i_1.$$
  Next, we choose the global phase reference and find absolute values of phases at every point in the network. Without loss of generality we can pick $\delta_2 = 0$, which gives $\delta_1 = \Delta_{1,2}$. Therefore, using the definition of the phase increment introduced earlier, $\delta = [\delta_0, \delta_1, \delta_3]^T$ can be found from the following system of chain equations:
  \begin{equation*}
  \begin{bmatrix}
  \Delta_{0,1} \\
  \Delta_{1,2} \\
  \Delta_{1,3}
  \end{bmatrix} = 
  \begin{bmatrix}
  1 & -1 & 0 \\
  0 & 1 & 0 \\
  0 & 1 & -1  
  \end{bmatrix}
  \begin{bmatrix}
  \delta_{0} \\
  \delta_{1} \\
  \delta_{3} 
  \end{bmatrix}.
  \end{equation*}
  The same calculation can be immediately generalised to the case when node $1$ has $n > 2$ links:
  \begin{equation*}
  \begin{bmatrix}
  \Delta_{0,1} \\
  \Delta_{1,2} \\
  \vdots \\
  \Delta_{1,n-1} \\
  \Delta_{1,n}
  \end{bmatrix} = 
  \begin{bmatrix}
  1 & -1 & 0 & \ldots & 0 & 0 \\
  0 & 1 & -1 & \ldots & 0 & 0 \\
  \vdots & \vdots & \vdots & \ddots & \vdots & \vdots \\
  0 & 1 & 0 & \ldots & -1 & 0 \\
  0 & 1 & 0 & \ldots & 0 & -1 \\
  \end{bmatrix}
  \begin{bmatrix}
  \delta_{0} \\
  \delta_{1} \\
  \vdots \\
  \delta_{n-1} \\
  \delta_{n}
  \end{bmatrix}.
  \end{equation*}
  It is known from graph theory, that for any tree graph the number of vertices is equal to the number of branches minus one \cite[Corollary $1.5.3$]{diestel2005graph}, i.e. for any tree graph $G(\mathcal{N}, \mathcal{E})$ with number of nodes/vertices $|\mathcal{N}|$ and 
  number of branches $|\mathcal{E}|$, $|\mathcal{N}| = |\mathcal{E}| - 1.$
  Thus any tree can be iteratively constructed by adding a pair consisting of a vertex and branch, which shows that the above procedure works for a chain feeder as well as for branching with $n$ links, therefore it is applicable for any tree network. 
\end{proof}
\begin{remark} \hfill
\begin{itemize}
  \item The phase matching procedure is required at branching nodes in a tree topology case, as briefly described in the proof above. For the complete algorithm refer to Algorithm \ref{phase_matching_procedure} in Section \ref{system_identification_section}.
  \item The forward model (\ref{lf-forward1}) is more cumbersome in a branching case. There are additional unknowns that need to be carried until termination. As an example, for the branching on Figure \ref{f2}:
  $$j^g_{0,1} = i_{1}e^{i\Delta_{0, 1}} + j_{1,2}^g + j_{1,3}^g,$$
  currents $j_{1,2}^g$ and $j_{1,3}^g$ are unknown, i.e. there are 2 equations, 4 unknowns.
\end{itemize}
\end{remark}

\subsection{Discussion}

As we saw in the previous sections, the forward model is not well suited to a complete low-voltage network description, hence we use the backward model further in the remainder of this paper. 

To conclude the modelling part we present an iterative model directly based on the backward model (\ref{l_pi_original}) for a single phase low-voltage network that we will exploit for impedance identification purposes. For any node $n$, the set of its ancestors ($\mathcal{A}_n$) and children ($\mathcal{C}_n$) in $G(\mathcal{N}, \mathcal{E})$ obey an iterative relationship:
\begin{equation}
\begin{split}
&{} v_{n}e^{i\Delta_{n, l}} - v_{l} = j_{n, l}z_{n, l}, \quad l \in \mathcal{C}_n,\\
& \sum_{k \in \mathcal{A}_n} j_{k, n} = i_{n} + \sum_{l \in \mathcal{C}_n} j_{n, l}e^{-i\Delta_{n, l}}.
\end{split}
\label{l_pi_general}
\end{equation}
Finally, it is worth noting that the system of equations (\ref{l_pi_general}) is not a new model for power networks. 
\begin{proposition}
  The system of equations (\ref{l_pi_general}) implies the standard power flow model:
  \begin{equation}
    \begin{split}
    &{} v_{n}^2 = v_{l}^2 + S_{n, l}z_{n, l}^* + S_{n, l}^*z_{n, l} + |j_{n, l}^gz_{n, l}|^2, \quad l \in \mathcal{C}_n, \\
    & \sum_{k \in \mathcal{A}_n} S_{k, n} - \sum_{l \in \mathcal{C}_n} (S_{n, l} + |j_{n, l}^g|^2z_{n, l})- s_{n} = 0.
    \end{split}
    \label{bpf_general}
  \end{equation}
  where $S_{n, l} = j_{n, l}^{g*}v_l^g$ and $s_{n} = i^{g*}_{n}v^{g}_{n}$.
\end{proposition}
\begin{proof}
  Exclude the phase term $e^{i\Delta_{n, l}}$ in (\ref{l_pi_general}) and introduce the \textit{receiving-end} \footnote{Note that the sign for terms with $|j_{n, l}|^2$ in resulting power flow equations is reversed with respect to the case when branch power is defined as sending-end complex power flow, i.e. $S_{n, l}=j_{n, l}^*v_n$.} complex power flow from $n$ to $l$ as $S_{n, l}=j_{n, l}^{g*}v_l^g = j_{n, l}^*|v_l|$ (since $ j_{n, l}^{g*}v_l^g = j_{n, l}^*e^{-i\delta_l}|v_l|e^{i\delta_l} = j_{n, l}^*|v_l|$), and  complex power injection $s_{n}$ to obtain:
  \begin{equation*}
    \begin{split}
    &{} v_{n}e^{i\Delta_{n, i}} - v_{i} = j_{n, i}z_{n, i}, \quad i \in \mathcal{C}_n,\\
    & \sum_{k \in \mathcal{A}_n} S_{k, n}^* = s_{n} + \sum_{i \in \mathcal{C}_n} (S_{n, i}^* + |j_{n, i}|^2z_{n, i}^*).
    \end{split}
  \end{equation*}
  Next, transfer the $v_i$ term to the right-hand side, take the magnitude squared of the first equation (also, note that $|j_{n, i}|^2 = |j_{n, i}^g|^2$) and complex conjugate of the second equation to obtain the well known power flow model (\ref{bpf_general}).
\end{proof} 

\begin{remark}
  Comparing equations (\ref{bpf_general}) and (\ref{l_pi_general}), we notice that the impedance term in the power flow model appears as $z_{n, l}$ and $|z_{n, l}|^2$, whereas the backward model has only $z_{n, l}$. Therefore, (\ref{l_pi_general}) is a better conditioned set of equations and we prefer it over (\ref{bpf_general}) for the problem of impedance identification, that we consider in detail in the next section.
\end{remark}

\section{System identification}
\label{system_identification_section}
In this section we consider the problem of system identification, that is, the problem of identifying impedances in a low-voltage network and we develop an approach that exploits the iterative structure of the backward model (\ref{l_pi_general}). We work through the details for a chain feeder model, however the proposed method can be generalised to any tree network as described earlier. We first formulate the optimisation problem for the whole network and then consider two methods to solve it approximately.

We start by considering the same chain feeder model as before (fig. \ref{f1}) with the following measurements available:
\begin{itemize}
  \item $\bm{v}_0$ - RMS \footnote{Root mean square} substation transformer voltage,
  \item $\bm{v}_{n} = \big[ v_{n}^{(1)}, \ldots, v_{n}^{(M)} \big]^T$ - RMS voltage measurements at node $n$,
  \item $|\bm{i}_{n}| = \big[ |i_{n}^{(1)}|, \ldots, |i_{n}^{(M)}| \big]^T$ - RMS consumption current measurements at node $n$, 
  \item $\bm{\theta}_{n} = \big[\theta_{n}^{(1)}, \ldots, \theta_{n}^{(M)}]^T$ - angle measurements between voltage and current phasors,
\end{itemize}
where indices $n = 1:N$ and $m = 1:M$ are used for node and measurement number respectively, e.g. $v_{n}^{(m)}$ - RMS value of the voltage at node $n$ corresponding to the $m$-th measurement\footnote{Although the index $m$ is used for measurement number, it does not necessarily mean that the corresponding quantity is measured, e.g. the value of $j_{n}^{(m)}$ corresponds to $m$-th measurement, but it is not measured directly and must be derived from other measurements.}. According to the minimum functionality requirements for smart metering infrastructure in Australia, every smart meter should be capable of measuring total active power, reactive power and voltage magnitude \cite[Table $6-2$]{SMmin_au}, which is equivalent to availability of $\bm{v}_{n}, |\bm{i}_{n}|, \bm{\theta}_n$ for all $n$. 

Regarding the measurements we assume the following.

\textbf{Assumption 1}. Smart meter measurements can be considered time synchronised with respect to system changes, i.e. for the $m$-th measurement taken from every node in the network, load impedances remain constant. 
This is justified because smart metering minimum functionality specifications require smart meters to maintain measurement time clocks to within a few seconds across a network, see e.g. \cite[Chapter $7.5$]{SMmin_au},\cite[Chapter $6$]{NMI_au}. 
\begin{remark}
Note that when PMU is available $\bm{v}^g_{n}$ and $\bm{i}^g_{n}$ values are measured directly. In contrast, the smart meter measurements ($\bm{v}_{n}, |\bm{i}_{n}|$) do not contain global phases $\bm{\delta}_n$. Nevertheless, the phase synchronism provided by PMU is not needed for our algorithms. 
\end{remark}
\textbf{Assumption 2} Line impedances and grid topology remain unchanged during the whole measurement process. 
\begin{remark} 
Algorithms that we develop in the next section allow us to monitor varying line impedances and, therefore, detect topology changes.
\end{remark}
\textbf{Assumption 3}. Smart meter measurement uncertainty is approximated by Gaussian noise $\mathcal{N}(0, \sigma)$, where $\sigma$ corresponds to $1\%$ full scale \footnote{i.e. the error value is taken with respect to the full scale of the measurement device} error which is practical for the majority of household smart meters, see e.g. technical reports \cite[Chapter $4$]{NMI_au}, \cite[Chapter $7.1$]{SMmin_au}.
\begin{remark}
  Note, however, that $1\%$ of full scale error for voltage, current and angle corresponds to about $2\%$ error in measured power, i.e. this setup can be considered as a worst case. Thus, we also consider $0.5\%$ and $0.1\%$ of full scale error cases.
\end{remark}

\subsection{Problem formulation. Chain feeder}

Consider a chain feeder (Fig. \ref{f1}) and formulate the problem of identifying impedances based on the backward model (\ref{l_pi_original}) written in multidimensional form (i.e. when multiple measurements are gathered):
\begin{equation}
\begin{split}
&{} \bm{v}_{n-1}e^{i\bm{\Delta}_{n}} - \bm{v}_{n} = \bm{j}_{n}z_{n}, \\
& \bm{j}_{n-1} = \bm{i}_{n-1} + \bm{j}_{n}e^{-i\bm{\Delta}_{n}},
\end{split}
\label{lf-sysid-model}
\end{equation}
where $\bm{v}_{n} = \big[ v^{(1)}_{n}, \ldots, v^{(M)}_{n} \big]^T$, $\bm{j}_{n} = \big[ j^{(1)}_{n}, \ldots, j^{(M)}_{n} \big]^T$, $\bm{\Delta_n} = \big[ \Delta^{(1)}_{n},$ $ \ldots, \Delta^{(M)}_{n} \big]^T$ and $e^{i\bm{\Delta}_{n}}$ is a component-wise operation as well as multiplication of two vectors $\bm{v}_{n-1}e^{i\bm{\Delta}_{n}}$ and $\bm{j}_{n}e^{i\bm{\Delta}_{n}}$.

Next, define a cost function for the $n$-th power line based on the first equation in (\ref{lf-sysid-model}): 
\begin{equation}
c_n(z_{n}, \bm{j}_n, \bm{\Delta}_{n}) = \Big\| \bm{v}_{n-1}e^{i\bm{\Delta}_{n}} - \bm{v}_{n} - \bm{j}_{n}z_{n} \Big\|_2^2,
\label{cost_function}
\end{equation} 
and introduce a cost-to-go function
\begin{equation}
J_{k:N} = \sum_{n = k}^{N} c_n(z_{n}, \bm{j}_n, \bm{\Delta}_{n})
\label{total_cost}
\end{equation}
for lines $N, N-1, \ldots, k$. Using the second equation in (\ref{lf-sysid-model}) we formulate the problem of impedance identification in a chain feeder:
\begin{equation}
\begin{aligned}
& \minimise_{z_{1:N}, \bm{\Delta}_{0:N-1}}
& & J_{1:N} \\
& \st
& & \bm{j}_{n-1} = \bm{i}_{n-1} + \bm{j}_{n}e^{-i\bm{\Delta}_{n}}, \\
&&& \bm{j}_{N} = \bm{i}_{N}, \\
&&& n = 1,\ldots, N-1.
\end{aligned}
\label{pf_network}
\end{equation}
The main advantage of this formulation is that the computation complexity scales linearly with the number of nodes, although the resulting optimisation problem is non-linear and non-convex. In the next subsections we consider two simplifications that allow solution for different practical scenarios. 

\subsection{Identification algorithm. Linearised case}

The first simplification is based on the following assumption which is very common in the analysis of low-voltage grids \cite{brice1982comparison}.

\textbf{Assumption 4} The phase increments $\Delta_n$ are negligibly small.  

\begin{remark}
Note, that we neglect the phase increments but we do not neglect the absolute phase at each node in the network. 
\end{remark}
Under Assumption 3, the model (\ref{lf-sysid-model}) takes the form:
\begin{equation}
  \begin{split}
  \bm{v}_{n-1} - \bm{v}_{n} = \bm{j}_{n}z_{n}, \\
  \bm{j}_{n-1} = \bm{i}_{n-1} + \bm{j}_{n}.
  \end{split}
  \label{lf-sysid-model-lin}
\end{equation}
Then the cost function (\ref{cost_function}) can be linearised and rewritten as follows by splitting it into real and imaginary parts:
\begin{equation}
  c^{l}_n(\bm{z}_{n}, \bm{j}_n) = \Big\| \bm{v}_{n-1} - \bm{v}_{n} - \bm{J}_{n}\bm{Q}_1\bm{z}_{n} \Big\|_2^2 + \Big\|\bm{J}_{n}\bm{Q}_2\bm{z}_{n} \Big\|_2^2,
  \label{cost_function-lin}
\end{equation}
where $\bm{J}_n \coloneqq \big[ \Real{\bm{j}_n} \quad \Imag{\bm{j}_n}\big]$, $\bm{z}_n \coloneqq \big[ \Real{z_n} \quad \Imag{z_n} \big]^T$, $\bm{Q}_1 \coloneqq 
\begin{bmatrix}
1 & 0 \\
0 & -1
\end{bmatrix}$ and $\bm{Q}_2 \coloneqq
\begin{bmatrix}
0 & 1 \\
1 & 0
\end{bmatrix}$.
Therefore, we obtain the optimisation problem of the form:
\begin{equation}
  \begin{aligned}
  & \minimise_{\bm{z}_{1:N}} 
  & & J^l_{1:N} \\
  & \st
  & & \bm{j}_{n-1} = \bm{i}_{n-1} + \bm{j}_{n}, \\
  &&& \bm{j}_{N} = \bm{i}_{N}, \\
  &&& n = 1,\ldots, N-1.
  \end{aligned}
  \label{lf-sysid-problem-lin}
\end{equation}
where $J^l_{1:N} = \sum_{n = 1}^{N} c^l_n(\bm{z}_{n}, \bm{j}_n)$ is the corresponding linearised cost-to-go function. 

This problem can be approached from the dynamic programming perspective.
% \begin{proposition}
% Minimising the sum in (\ref{lf-sysid-problem-lin}) is equivalent to minimising each term $c^l_n(\bm{z}_{n}, \bm{j}_n)$ separately.
% \end{proposition}
% \begin{proof}
  To show this, consider the cost-to-go function for the last two power lines:
  \begin{equation*}
    \begin{split}
    &{}J_{N-1:N} = \min_{\bm{z}_{N-1:N}} \big[ c^l_{N-1}(\bm{z}_{N-1}, \bm{j}_{N-1}) + c^l_N(\bm{z}_{N}, \bm{j}_{N}) \big] = \\
    & = \min_{\bm{z}_{N-1}} \big[ c_{N-1}(\bm{z}_{N-1}, \bm{i_{N-1}} + \bm{j}_{N}) \big] + \min_{\bm{z}_{N}} \big[ c^l_N(\bm{z}_{N}, \bm{j}_{N}) \big],
    \end{split}
  \end{equation*}
  where both terms are the impedance identification problems for lines $N$ and $N-1$. Using the same argument for all power lines in the feeder, we conclude that minimising sum in (\ref{pf_network}) is equivalent to minimising each term separately according to Algorithm \ref{a2}. Since it identifies the network starting from the edge node and propagating all the way through to the first node, we call it LBCI (Linearised Backward Calculation of Impedances).
% \end{proof}
\begin{algorithm}[H]
\caption{LBCI algorithm for a chain feeder}
\label{a2}
\begin{algorithmic}[1]
\State $\mathcal{M} \coloneqq \{\bm{v}_{n}, |\bm{i}_{n}|, \bm{\theta}_{n} ~\big|~ n=1, \ldots, N \} \cup \{ \bm{v}_0 \}$
\Function{LBCI}{$\mathcal{M}$}
\State $\bm{j}_{N} = \bm{i}_{N}$
\For{$n = N : 1$}
\State $\widehat{\bm{z}}_n = \argmin_{\bm{z}_{n}} c^{l}_n(\bm{z}_{n}, \bm{j}_n)$
\State $\bm{j}_{n-1} = \bm{i}_{n-1} + \bm{j}_{n}$
\EndFor
\State \Return $\widehat{\bm{Z}} = \big[ \widehat{\bm{z}}_{1:N} \big], \bm{J} = \big[ \bm{j}_{1:N} \big]$
\EndFunction
\end{algorithmic}
\end{algorithm} 

\begin{remark}
\hfill
\begin{itemize}
  \item Minimisation of $c^l_n(\bm{z}_n, \bm{j}_n)$ with respect to $\bm{z}_n$ is a simple unconstrained least squares problem, where the second term can be interpreted as a regularisation term. Note, that it appears naturally in the problem formulation.
  \item Formally, measurement noise affects the $\bm{J}_n$ term, since it contains data (real and imaginary parts of $\bm{j}_n$) derived from measurements. In practice, the main source of error is the voltage measurement noise, therefore the error in matrix $\bm{J}_n$ can be neglected in favour of using a linear least squares approach. Linear least squares estimation is also preferable as a simple (from a computational perspective for decentralised implementation) and flexible (i.e. has recursive, weighted modifications) method. This allows real-time tracking of impedances in the network as new measurements become available. Nevertheless, when errors in $\bm{J}_n$ need to be taken into account, total least squares (TLS) techniques can be used \cite{golub1999tikhonov}.
  \item There is widely used variation of LBCI algorithm where the regularisation part is ignored in (\ref{cost_function-lin}), i.e.:
  \begin{equation*}
    c^{l}_n(\bm{z}_{n}, \bm{j}_n) = \Big\| \bm{v}_{n-1} - \bm{v}_{n} - \bm{J}_{n}\bm{Q}_1\bm{z}_{n} \Big\|_2^2
  \end{equation*}
  We refer to this algorithm as LBCI-old. In \cite{peppanen2015distribution} a very similar approach is used, although the main difference from the case we consider is computation complexity. The authors propose solving a least-squares problem with measurement matrix of size $MN \times (M+2N)$ which has complexity around $O((M+2N)^2 MN)$ whereas in our case we solve $N$ problems with $M \times 2$ matrices which results in $O(MN)$.
\end{itemize}
\end{remark}

Generalisation to a tree topology feeder involves one extra step: we need to keep track of phase increments in order to match the absolute phases of line currents at the intersection point. Note, that because of linearisation we do not obtain correct values of absolute phases, however this is not an issue for identification purposes, since the backward model (\ref{l_pi_original}) is based on the incremental phase information. 

Let $\mathcal{M}^{i}$ and $\mathcal{M}^{j}$ be the sets of measurements for branches $i$ and $j$ that have a common node $n^*$ and sets of indices $\mathcal{N}^i$ and $\mathcal{N}^j$ correspondingly. We briefly summarise the procedure of phase matching:
\begin{algorithm}[H]
\caption{Phase matching procedure}
\label{phase_matching_procedure}
\begin{algorithmic}[1]
%\Require $\mathcal{M}^i$, $\mathcal{M}^j$
\State Calculate: $$\widehat{\bm{Z}}^i, \bm{J}^i =\textproc{LBCI}(\mathcal{M}^i), \quad \widehat{\bm{Z}}^j, \bm{J}^j = \textproc{LBCI}(\mathcal{M}^j).$$
\State For a common node $n^*$ find absolute phases as a sum of phase increments $\bm{\Delta}_{k}$ \footnotemark:
$$\bm{\Delta}_{k} = \frac{\bm{J}_{k}\bm{Q}_2 \widehat{\bm{z}}^i_{k}}{\bm{v}_{k-1}}, \quad \bm{\delta}^i_{n*} \approx \sum_{k \in \mathcal{N}_i} \bm{\Delta}_{k} , \quad \bm{\delta}^j_{n*} \approx \sum_{k \in \mathcal{N}_j} \bm{\Delta}_{k}.$${}
\State Chose $\bm{\delta}^i_{n*}$ as a reference and update phases in branch $j$ by a factor of $\bm{\delta}^i_{n*} - \bm{\delta}^j_{n*}$. This is shown by black arrows in Fig. \ref{f4}.
%\Ensure $\widehat{\bm{Z}} = \big[ \widehat{\bm{z}}_{1:N} \big], \bm{J} = \big[ \bm{j}_{1:N} \big]$
\end{algorithmic}
\end{algorithm}
\footnotetext{Formula for $\bm{\Delta}_{k}$ is obtained from (\ref{lf-sysid-model}) after taking imaginary part of the first equation.}
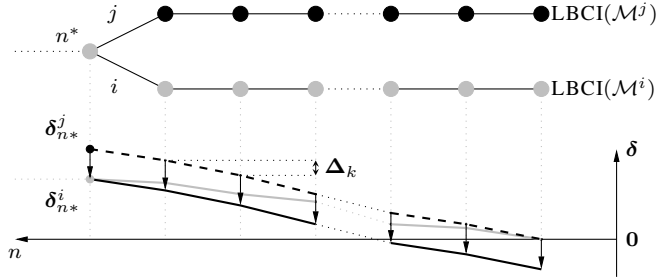
\begin{figure}[ht]
  \centering
  \begin{tikzpicture}
    \coordinate (a0) at (-1, 3);
    \coordinate (a1) at (1, 2.5); 
    \coordinate (a2) at (2, 2.5);
    \coordinate (a3) at (3, 2.5); 
    \coordinate (a4) at (4, 2.5); 
    \coordinate (a5) at (5, 2.5); 
    \coordinate (a6) at (6, 2.5); 

    \coordinate (b0) at (0, 3); 
    \coordinate (b1) at (1, 3.5); \draw[fill] (b1) circle [radius=0.1];
    \coordinate (b2) at (2, 3.5); \draw[fill] (b2) circle [radius=0.1];
    \coordinate (b3) at (3, 3.5); \draw[fill] (b3) circle [radius=0.1];
    \coordinate (b4) at (4, 3.5); \draw[fill] (b4) circle [radius=0.1];
    \coordinate (b5) at (5, 3.5); \draw[fill] (b5) circle [radius=0.1];
    \coordinate (b6) at (6, 3.5); \draw[fill] (b6) circle [radius=0.1];

    \coordinate (a0_) at (0, 1.3); \draw[fill] (a0_) circle [radius=0.01];
    \coordinate (a1_) at (1, 1.25); \draw[fill] (a1_) circle [radius=0.01];
    \coordinate (a2_) at (2, 1.1); \draw[fill] (a2_) circle [radius=0.01];
    \coordinate (a3_) at (3, 1); \draw[fill] (a3_) circle [radius=0.01];
    \coordinate (a4_) at (4, 0.7); \draw[fill] (a4_) circle [radius=0.01];
    \coordinate (a5_) at (5, 0.65); \draw[fill] (a5_) circle [radius=0.01];
    \coordinate (a6_) at (6, 0.5); \draw[fill] (a6_) circle [radius=0.01];

    \coordinate (b0_) at (0, 1.7); \draw[fill] (b0_) circle [radius=0.01];
    \coordinate (b1_) at (1, 1.55); \draw[fill] (b1_) circle [radius=0.01];
    \coordinate (b2_) at (2, 1.35); \draw[fill] (b2_) circle [radius=0.01];
    \coordinate (b3_) at (3, 1.1); \draw[fill] (b3_) circle [radius=0.01];
    \coordinate (b4_) at (4, 0.85); \draw[fill] (b4_) circle [radius=0.01];
    \coordinate (b5_) at (5, 0.7); \draw[fill] (b5_) circle [radius=0.01];
    \coordinate (b6_) at (6, 0.5); \draw[fill] (b6_) circle [radius=0.01];

    \coordinate (y0) at (7, 0);
    \coordinate (y1) at (7, 1.7);
    \coordinate (x0) at (7, 0.5);
    \coordinate (x1) at (-1, 0.5);

    \draw[dotted] (a0) -- (b0) node[above left]{\footnotesize $n^*$};
    \draw (b0) -- node[above left]{\footnotesize $j$}(b1) -- (b2) -- (b3);
    \draw (b0) -- node[below left]{\footnotesize $i$}(a1) -- (a2) -- (a3);
    \draw[dotted] (b3) -- (b4);
    \draw[dotted] (a3) -- (a4);
    \draw (b4) -- (b5) -- (b6) node[right]{\footnotesize $\textproc{LBCI}(\mathcal{M}^j)$};
    \draw (a4) -- (a5) -- (a6) node[right]{\footnotesize $\textproc{LBCI}(\mathcal{M}^i)$};

    \draw[gray!50, dotted] (b0) -- (0, 0.5);
    \draw[gray!50, dotted] (a1) -- (1, 0.5);
    \draw[gray!50, dotted] (a2) -- (2, 0.5);
    \draw[gray!50, dotted] (a3) -- (3, 0.5);
    \draw[gray!50, dotted] (a4) -- (4, 0.5);
    \draw[gray!50, dotted] (a5) -- (5, 0.5);
    \draw[gray!50, dotted] (a6) -- (6, 0.5);

    \draw[gray!50, fill] (b0) circle [radius=0.1];
    \draw[fill, gray!50] (a1) circle [radius=0.1];
    \draw[fill, gray!50] (a2) circle [radius=0.1];
    \draw[fill, gray!50] (a3) circle [radius=0.1];
    \draw[fill, gray!50] (a4) circle [radius=0.1];
    \draw[fill, gray!50] (a5) circle [radius=0.1];
    \draw[fill, gray!50] (a6) circle [radius=0.1];
    
    \draw[arrows = {-Stealth[inset=0pt, angle=25:5pt]}] (y0) -- (y1) node[right]{\footnotesize $\bm{\delta}$}; \draw[arrows = {-Stealth[inset=0pt, angle=25:5pt]}] (x0) node[right]{\footnotesize $\bm{0}$} -- (x1)node[below]{\footnotesize $n$};

    \draw[gray!50, thick] (a0_) -- (a1_) -- (a2_) -- (a3_); \draw[gray!50, dotted] (a3_) -- (a4_); \draw[gray!50, thick] (a4_) -- (a5_) -- (a6_);
    \draw[gray!50, dotted] ($(a0_) - (1, 0)$) -- (a0_);

    \draw[thick, dashed] (b0_) node[above left]{\footnotesize $\bm{\delta}^j_{n*}$} -- (b1_) -- (b2_) -- (b3_); \draw[dotted] (b3_) -- (b4_); \draw[thick, dashed] (b4_) -- (b5_) -- (b6_);
    \draw[fill] (b0_) circle [radius=0.05];
    \draw[dotted] (b1_) -- ($(b1_) + (2, 0)$);
    \draw[dotted] (b2_) -- ($(b2_) + (1, 0)$);
    \draw[<->, >=stealth, line width=0.1pt] ($(b1_) + (2, 0)$) -- node[right]{\footnotesize $\bm{\Delta}_{k}$} ($(b2_) + (1, 0)$);

    \coordinate (shift) at ($(a0_) - (b0_)$); 
    \draw[thick] ($(b0_) + (shift)$) node[below left]{\footnotesize $\bm{\delta}^i_{n*}$}-- ($(b1_) + (shift)$) -- ($(b2_) + (shift)$) -- ($(b3_) + (shift)$); \draw[dotted] ($(b3_) + (shift)$) -- ($(b4_) + (shift)$); \draw[thick] ($(b4_) + (shift)$) -- ($(b5_) + (shift)$) -- ($(b6_) + (shift)$);
    \draw[fill, gray!50] ($(b0_) + (shift)$) circle [radius=0.05];

    \draw[arrows = {-Stealth[inset=0pt, angle=25:5pt]}] (b0_) -- ($(b0_) + (shift)$);
    \draw[arrows = {-Stealth[inset=0pt, angle=25:5pt]}] (b1_) -- ($(b1_) + (shift)$);
    \draw[arrows = {-Stealth[inset=0pt, angle=25:5pt]}] (b2_) -- ($(b2_) + (shift)$);
    \draw[arrows = {-Stealth[inset=0pt, angle=25:5pt]}] (b3_) -- ($(b3_) + (shift)$);
    \draw[arrows = {-Stealth[inset=0pt, angle=25:5pt]}] (b4_) -- ($(b4_) + (shift)$);
    \draw[arrows = {-Stealth[inset=0pt, angle=25:5pt]}] (b5_) -- ($(b5_) + (shift)$);
    \draw[arrows = {-Stealth[inset=0pt, angle=25:5pt]}] (b6_) -- ($(b6_) + (shift)$);

  \end{tikzpicture}
  \caption{Phase matching procedure example. Top: two branches $i$ (grey dots) and $j$ (black dots). Bottom: corresponding phase plots. The grey line is the phase of the reference branch $i$, dotted line shows the initial phase for branch $j$, black solid line is the updated phase for branch $j$.}
  \label{f4}
\end{figure}

\noindent The illustration of this procedure is given in Figure \ref{f4}, where the black dotted line, corresponding to the phases along $j$-th branch before matching, gets shifted by the amount of phase difference at the node $n^*$. The complete algorithm for a tree topology case can be obtained by combining this procedure with the chain feeder algorithm. Therefore, under the linearisation assumption (Assumption 3), the problem of impedance identification can be solved optimally. 

To conclude, the developed algorithm is fast, simple and can be implemented in a distributed fashion, as discussed in Section \ref{practical_modifications_section}. However, we lose information about phase angle increments at each power line which leads to unknown absolute phase at each node. It turns out that we can address this issue very efficiently without losing the key properties of the developed algorithm. 

\subsection{Identification algorithm. Non-linear case}

In this section we consider an alternative simplification of (\ref{pf_network}) where minimisation of each cost function (\ref{cost_function}) is a non-linear and non-convex problem. Noticeably, any such problem has certain properties that can be exploited in order to solve it globally. Although the modification we introduce does not lead to a global solution of the original non-linear problem (\ref{pf_network}), we show that any solution it finds is closer to the global optimum of (\ref{pf_network}) than that found in linearised case.

Consider the cost function (\ref{cost_function}) and, using the notation introduced in (\ref{cost_function-lin}), write it as follows by splitting real and imaginary parts:
\begin{equation}
\begin{split}
c_n(\bm{z}_{n}, \bm{j}_n, \bm{\gamma}_{n}) ={}& \Big\| \bm{V}_{n-1}\bm{\gamma}_{n} - \bm{v}_{n} - \bm{J}_{n}\bm{Q}_1\bm{z}_{n} \Big\|_2^2 + \\
& \Big\|\bm{V}_{n-1}\sqrt{\bm{1} - \bm{\gamma}^2_{n}} - \bm{J}_{n}\bm{Q}_2\bm{z}_{n} \Big\|_2^2,
\end{split}
\label{cost_function-nlin}
\end{equation}
where $\bm{\gamma}_n \coloneqq \cos{\bm{\Delta}_{n}}$, $\bm{V}_n \coloneqq \diag{(\bm{v}_{n})}$, $\bm{1} = \big[1, \ldots, 1 \big]^T$, operations $\sqrt{(\cdot)}$, $\cos{(\cdot)}$ and $(\cdot)^2$ are component-wise. 

Next, instead of minimising the original cost $c_n(\bm{z}_{n}, \bm{j}_n, \bm{\gamma}_{n})$ at each step, we minimise the first term and use the second term as a constraint:
\begin{equation}
\begin{aligned}
& \minimise_{\bm{\gamma}_{n}, \bm{z}_{n}} 
& & \Big\| \bm{V}_{n-1}\bm{\gamma}_{n} - \bm{v}_{n} - \bm{J}_{n}\bm{Q}_1\bm{z}_{n} \Big\|_2^2 \\
& \st
& & \bm{\gamma}_{n} = \sqrt{\bm{1} - \frac{(\bm{J}_{n}\bm{Q}_2\bm{z}_{n})^2}{\bm{v}^2_{n-1}}},
\end{aligned}
\label{sysid-nlin}
\end{equation}
where operations $\frac{(\cdot)}{(\cdot)}$ and $(\cdot)^2$ are component-wise. The reasoning behind this formulation will become clear in the next subsection where we make comparison with the LBCI algorithm.
\begin{remark}
In the constraint equation for $\bm{\gamma}_{n}$ we use only the positive root. By doing so, we restrict consideration for phase angles $\Delta^{(m)}_{n} \in [-\pi/2; \pi/2]$ (or $0 \le \gamma^{(m)}_{n} \le 1$) which is practical for low-voltage power networks.
\end{remark}

The main idea for the optimisation algorithm is based on Lemma \ref{lemma1} (see Appendix). In the context of problem (\ref{sysid-nlin}), let:
\begin{equation*}
\begin{split}
{}& f(\bm{\gamma}_{n}, \bm{z}_n) \coloneqq \Big\| \bm{V}_{n-1}\bm{\gamma}_{n} - \bm{v}_{n} - \bm{J}_{n}\bm{Q}_1\bm{z}_{n} \Big\|_2^2, \\
& \bm{g}(\bm{z}_{n}) \coloneqq \sqrt{\bm{1} - \frac{(\bm{J}_{n}\bm{Q}_2\bm{z}_{n})^2}{\bm{v}^2_{n-1}}}, \\
& \bm{h}(\bm{\gamma}_{n}) \coloneqq \big[ \bm{J}_{n}\bm{Q}_1 \big]^{\dagger} (\bm{V}_{n-1}\bm{\gamma}_{n} - \bm{v}_{n}).
\end{split}
\end{equation*}
First, we check if the requirements of Lemma \ref{lemma1} are satisfied:
\begin{itemize}
\item $\bm{V}_{n-1}$ is a full rank matrix since $\bm{V}_{n-1} = \diag{(\bm{v}_{n-1})}$. Regarding $\bm{J}_{n}$, we assume that it is a full rank matrix, in other words, there should be enough variation in measurements.
\item $\bm{g}(\bm{z}_{n})$ is a surjective mapping for all $\bm{z}_{n}$ that satisfy:
$$(\bm{J}_{n}\bm{Q}_2\bm{z}_{n})^2 \le \bm{v}^2_{n-1}$$
This condition requires that the voltage drop at the line is less than or equal to the voltage at its node, which is always true for low-voltage networks.
\end{itemize}
\begin{remark}
  $\bm{J}_n$ is not a full rank matrix if it contains identical measurements or if the line current is zero for all measurements except one. Both cases are pathological, and to handle them in practise, we introduce regularisation term into the problem formulation (see the next subsection). There is another case when $\bm{J}_n$ becomes singular: a power factor among all loads in the feeder is the same and remains constant, which is also highly improbable. We address this problem in simulations (see section \ref{num_sim_section}). 
  \label{remark_nonzero_J}
\end{remark}

Second, according to Lemma \ref{lemma1} we need to find a fixed point of $\bm{g} \circ \bm{h}$. 
\begin{proposition}
A fixed point of $\bm{g} \circ \bm{h}$ can be found via the following iteration:
\begin{equation}
\begin{split} 
\bm{\gamma}_{n}^{(i+1)} = \bm{\gamma}_{n}^{(i)} + \alpha\Big[\bm{g} \circ \bm{h}( \bm{\gamma}_{n}^{(i)} ) - \bm{\gamma}_{n}^{(i)}\Big], ~\bm{\gamma}_{n}^{(0)} = \bm{1}, 
\end{split}
\label{iteration}
\end{equation}
where $0 < \alpha < 1$.
\end{proposition}
\begin{proof}
  The iteration (\ref{iteration}) is a gradient descent type algorithm for the area maximisation problem. To see this, consider:
  \begin{equation*}
    \begin{split}
    {}& \nabla_{\bm{\gamma}_{n}} \Big[ \int_{\bm{\gamma}_{n}^{(0)}}^{\bm{\gamma}_{n}} \big[\bm{g} \circ \bm{h}(\bm{\bar{\gamma}}_{n}) - \bm{\bar{\gamma}}_{n}\big] d \bm{\bar{\gamma}_{n}} \Big] = \bm{g} \circ \bm{h}(\bm{\gamma}_{n}) - \bm{\gamma}_{n}.
    \end{split}
  \end{equation*}
  As follows from the properties of definite integrals, the function in brackets is strictly increasing and therefore standard optimisation techniques can be applied in order to find its global maximum. Figure \ref{fig1} illustrates this idea for some scalar function $g \circ h$. 
  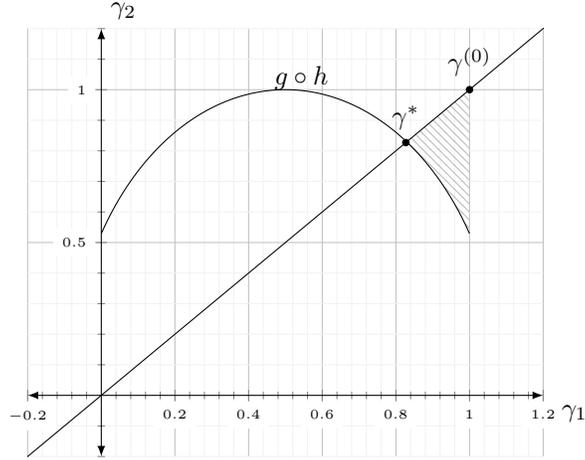
\begin{figure}[t]
  \centering
  \begin{tikzpicture}
  \begin{axis}[xmin=0,xmax=1,
              ymin=0,ymax=1,
              grid=both,
              grid style={line width=.1pt, draw=gray!10},
              major grid style={line width=.2pt,draw=gray!50},
              axis lines=middle,
              minor tick num=4,
              enlargelimits={abs=0.2},
              axis line style={latex-latex},
              ticklabel style={font=\tiny, fill=white},
              xlabel style={at={(ticklabel* cs:1)},anchor=north west},
              ylabel style={at={(ticklabel* cs:1)},anchor=south west},
              xlabel = {$~\gamma_1$},
              ylabel = {$\gamma_2$}]
  \addplot[name path=F,black!,domain={0:1}, samples=101] {sqrt(1 - 0.02*(12*x - 6)^2)} node[pos=0.6, above left]{$g \circ h$};
  %Sqrt[1 - 0.0000135558 (-200 + 220 x)^2] mathematica approximate result
  \addplot[name path=G,black,domain={-1:2}] {x};
  \addplot[pattern=north west lines, pattern color=gray!50]fill between[of=F and G, soft clip={domain=0.82:1}];
  \node[label={90:{$\gamma^{(0)}$}},circle,fill,inner sep=1pt] at (axis cs:1,1) {};
  \node[label={90:{$\gamma^*$}},circle,fill,inner sep=1pt] at (axis cs:0.827,0.827) {};
  \end{axis}
  \end{tikzpicture}
  \caption{Area maximisation problem for a scalar case. The objective is to maximise the area between $g \circ h$ and a line, starting from $\gamma^{(0)}$. The result is shown in dashed gray.}
  \label{fig1}
  \end{figure}

  % Rewrite (\ref{iteration}) in the form of the Richardson iteration \cite{kelleyiterative}:
  % \begin{equation*}
  %   \bm{\gamma}_{n}^{(i+1)} = \bm{\gamma}_{n}^{(i)} (\bm{I} - \bm{A}) + \alpha \bm{g} \circ \bm{h}( \bm{\gamma}_{n}^{(i)})
  % \end{equation*}
  % where $\bm{A} \coloneqq \alpha \bm{I}$ and $\bm{I}$ is the identitiy matrix. 
  % It follows that $\|\bm{I} - \bm{A}\|_2 < 1$ since $0 < \alpha < 1$. Therefore applying Corollary 1.2.1 of the Lemma 1.2.1 from \cite{kelleyiterative} the iteration (\ref{iteration}) converges to a fixed point $\bm{\gamma}_{n}^{*} = \bm{g} \circ \bm{h}( \bm{\gamma}_{n}^{*})$ for any initial $\bm{\gamma}_{n}^{(0)}$.
\end{proof}
\begin{remark} \hfill
\begin{itemize}
  \item The sign of $\alpha$ defines the direction of integration. For instance, in (\ref{iteration}) the direction is chosen so, that the iteration converges to the $\bm{0} \le \bm{\gamma}_{n}^{*} \le \bm{\gamma}_{n}^{(0)}$. Importantly, $\bm{\gamma}_{n}^{*}$ that solves (\ref{sysid-nlin}) is always non-negative, this property allows using iteration (\ref{iteration}) and obtain a global solution. 
  \item In the context of power systems, at the first step ($i = 0$) of the iteration (\ref{iteration}) we solve the minimisation problem (\ref{sysid-nlin}) with chain approximation around $\bm{\gamma}_n = \bm{1}$. This is equivalent to Assumption 2(see previous subsection). At the next iteration the algorithm consequently improves previous estimation converging to the solution. 
\end{itemize}
\end{remark}

Finally, we formalise the resulting identification algorithm with iteration (\ref{iteration}) for a single power line between nodes $n-1$ and $n$ in Algorithm \ref{a3}. 
\begin{algorithm}[ht]
\caption{Local iteration for a single power line.
}
\label{a3}
\begin{algorithmic}[1]
\State $\mathcal{M}_{n-1:n} \coloneqq \{\bm{v}_{k}, |\bm{i}_{k}|, \bm{\theta}_{k} ~\big|~ k=n, n-1 \}$
\Function{LineCalc}{$\mathcal{M}_{n-1:n}, \epsilon, \alpha, \bm{j}_n$}
%\Comment Data, accuracy and step size
\State $i = 0$, $\bm{\gamma}_{n}^{(0)} = \mathbf{1}$
\Do 
\State $\widehat{\bm{z}}_n^{(i)} = \bm{h}(\bm{\gamma}_{n}^{(i)})$
\State $\bm{\gamma}_{n}^{(i+1)} = \bm{\gamma}_{n}^{(i)} - \alpha\Big[\bm{g}(\widehat{\bm{z}}_n^{(i)}) - \bm{\gamma}_{n}^{(i)}\Big]$
%\State $i = i+1$
\doWhile{$\|\bm{g}(\widehat{\bm{z}}_n^{(i)}) - \bm{\gamma}_{n}^{(i)}\|_2 > \epsilon$}
\State \Return $\widehat{z}, \bm{\gamma}$
\EndFunction
\end{algorithmic}
\end{algorithm} 
Note, that this algorithm uses information from the previous line of the network, i.e. it needs $\bm{j}_n$ to be known. Therefore, using the backward model (\ref{l_pi_original}) iteratively we formulate the complete algorithm for a chain feeder identification in Algorithm \ref{a4}
\begin{algorithm}[H]
\caption{BCI algorithm for a chain feeder.}
\label{a4}
\begin{algorithmic}[1]
\State $\mathcal{M} \coloneqq \{\bm{v}_{n}, |\bm{i}_{n}|, \bm{\theta}_{n} ~\big|~ n=1, \ldots, N \} \cup \{ \bm{v}_0 \}$
\Function{BCI}{$\mathcal{M}, \epsilon, \alpha$}
\State $\bm{j}_{N} = \bm{i}_{N}$
\For{$n = N : 1$}
\State $\widehat{\bm{z}}_n, \bm{\gamma}_n = $ \textproc{LineCalc}($\mathcal{M}_{n-1:n}, \epsilon, \alpha, \bm{j}_n$)
\State $\bm{j}_{n-1} = \bm{i}_{n-1} + \bm{j}_{n}\frac{|v_{n}|}{|v_{n-1}|} + \widehat{z}_n^*\frac{|j_{n}|^2}{|v_{n-1}|}$ 
\Comment $(\ast)$ 
\EndFor
\State \Return $\bm{\gamma}, \widehat{\bm{z}}, \bm{j}$
\EndFunction
\end{algorithmic}
\end{algorithm} 
\begin{remark}
Equation $(\ast)$ in the Algorithm \ref{a4} can be obtained from the backward model (\ref{l_pi_original}) after excluding term $e^{i\Delta_n}$. 
\end{remark}

\subsection{Discussion}

The main objective of impedance identification algorithm is to minimise a cost function given by (\ref{cost_function-nlin}) for each power line. In this section we compare the optimality of solutions found by BCI and LBCI algorithms.

Again, consider a power line connecting nodes $n-1$ and $n$. The following proposition holds for solutions found by LBCI and BCI: 
\begin{proposition}
\label{prop_bci_vs_lbci}
Let $(\bm{z}_n^{LBCI}, \bm{1})$ is the solution found by the LBCI algorithm and $(\bm{z}_n^{BCI}, \bm{\gamma}^{BCI}_{n})$ is the solution found by the BCI algorithm.
Also, let $N(\bm{z}_n, \bm{\gamma}_{n}) \coloneqq \big\| \bm{V}_{n-1}\bm{\gamma}_{n} - \bm{v}_{n} - \bm{J}_{n}\bm{Q}_1\bm{z}_{n} \big\|_2^2 + \big\|\bm{V}_{n-1}\sqrt{\bm{1} - \bm{\gamma}^2_{n}} - \bm{J}_{n}\bm{Q}_2\bm{z}_{n} \big\|_2^2$ then
$$N(\bm{z}_n^{BCI}, \bm{\gamma}^{BCI}_{n}) \le N(\bm{z}_n^{LBCI}, \bm{1}),$$
i.e. the BCI algorithm finds the solution that is closer to the optimum of (\ref{pf_network}) than that of LBCI.
\end{proposition}
\begin{proof}
See appendix.
\begin{comment}
Consider a linearised BCI algorithm for the power line between nodes $n-1$ and $n$. The solution it finds satisfies:
\begin{equation*}
\begin{split}
{}& N(\bm{z}_n^{LBCI}, \bm{1}) = \big\| \bm{V}_{n-1} - \bm{v}_{n} - \bm{J}_{n}\bm{Q}_1\bm{z}^{LBCI}_{n} \big\|_2^2 + \\
& + \big\|\bm{J}_{n}\bm{Q}_2\bm{z}^{LBCI}_{n} \big\|_2^2 \ge \big\| \bm{V}_{n-1} - \bm{v}_{n} - \bm{J}_{n}\bm{Q}_1\bm{z}^{LBCI}_{n} \big\|_2^2.
\end{split}
\end{equation*}

The corresponding expressions for $\bm{z}_n^{BCI}$ and $\bm{z}_n^{LBCI}$ are given by:  
\begin{equation*}
\begin{split}
{}& \bm{z}_n^{BCI} = \big[ \bm{J}_{n}\bm{Q}_1 \big]^{\dagger} (\bm{V}_{n-1}\bm{\gamma}^{BCI}_{n} - \bm{v}_{n}), \\
& \bm{z}_n^{LBCI} = \big[ \bm{J}_{n}\bm{Q}_1 \big]^{\dagger} (\bm{V}_{n-1} - \bm{v}_{n}),
\end{split}
\end{equation*}
where $\bm{\gamma}^{BCI}_{n} \le \bm{1}$. Therefore:
\begin{equation*}
\begin{split}
{}& N(\bm{z}_n^{BCI}, \bm{\gamma}^{BCI}_{n}) = \big\| \bm{V}_{n-1}\bm{\gamma}^{BCI}_{n} - \bm{v}_{n} - \bm{J}_{n}\bm{Q}_1\bm{z}^{BCI}_{n} \big\|_2^2 \le\\
& \Big\|\big[ \bm{I} - \bm{J}_{n}\bm{Q}_1\big[ \bm{J}_{n}\bm{Q}_1 \big]^{\dagger}\big](\bm{V}_{n-1}\bm{\gamma}^{BCI}_{n} - \bm{v}_{n})\Big\|_2^2 \le \\
& \Big\|\big[ \bm{I} - \bm{J}_{n}\bm{Q}_1\big[ \bm{J}_{n}\bm{Q}_1 \big]^{\dagger}\big](\bm{V}_{n-1} - \bm{v}_{n})\Big\|_2^2 = \\
& \big\| \bm{V}_{n-1} - \bm{v}_{n} - \bm{J}_{n}\bm{Q}_1\bm{z}^{LBCI}_{n} \big\|_2^2 \le N(\bm{z}_n^{LBCI}, \bm{1}).
\end{split}
\end{equation*}
\end{comment}
\end{proof}
\begin{corollary} \label{corollary} In the presence of high measurement noise, we can introduce a regularisation term (as in the LBCI case) for BCI algorithm so that:
$N(\bm{z}_n^{BCI}, \bm{\gamma}^{BCI}_{n}) \le N(\bm{z}_n^{LBCI}, \bm{1})$. 
\begin{proof}
See appendix.
\end{proof}

\begin{comment}
To see this, note that:
\begin{equation*}
\begin{split}
{}& N(\bm{z}_n^{BCI}, \bm{\gamma}^{BCI}_{n}) \le \Big\|\big[ \bm{I} - \bm{J}_{n}\bm{Q}_1\big[ \bm{J}_{n}\bm{Q}_1 \big]^{\dagger}\big](\bm{V}_{n-1} - \bm{v}_{n})\Big\|_2^2 \le \\
& \Big\|\big[ \bm{I} - \bm{J}_{n}\bm{Q}_1\big[ \bm{J}_{n}\bm{Q}_1 \big]^{\dagger}\big](\bm{V}_{n-1} - \bm{v}_{n})\Big\|_2^2 + \mu \big\|\bm{J}_{n}\bm{Q}_2\bm{z}^{BCI}_{n} \big\|_2^2 \le \\
& \Big\|\big[ \bm{I} - \bm{J}_{n}\bm{Q}_1\big[ \bm{J}_{n}\bm{Q}_1 \big]^{\dagger}\big](\bm{V}_{n-1} - \bm{v}_{n})\Big\|_2^2 + \big\|\bm{J}_{n}\bm{Q}_2\bm{z}^{LBCI}_{n} \big\|_2^2 = \\
& N(\bm{z}_n^{LBCI}, \bm{1}),
\end{split}
\end{equation*}
where $0 \le \mu \le 1$ is a regularisation parameter.
\end{comment}
\end{corollary}
Importantly, all the key properties of the LBCI algorithm are inherited by the BCI algorithm, since $\bm{h}(\bm{\gamma}_{n}^{(i)})$ is a least squares solution of the unconstrained problem (\ref{sysid-nlin}) for a given $\bm{\gamma}_{n}^{(i)}$. 

\subsection{Practical modifications}
\label{practical_modifications_section}
This subsection describes some additional modifications that improve the LBCI/BCI algorithm performance. 

We start from the idea of decentralised calculations, whereby the BCI algorithm above can be executed by smart meters in a decentralised manner with a modification to take advantage of a known $X/R$ ratio\footnote{$X/R$ ratio is the ratio between inductance and resistance of the cable} from cable data. 

\textbf{Regularisation}. As was mentioned in the section \ref{remark_nonzero_J}, if elements in $\bm{J}_n$ are equal to zero (or close to it) then $\bm{J}_n$ is no longer a full rank. To address this issue we can add regularisation term to the objective function in (\ref{sysid-nlin}), i.e:
$$\Big\| \bm{V}_{n-1}\bm{\gamma}_{n} - \bm{v}_{n} - \bm{J}_{n}\bm{Q}_1\bm{z}_{n} \Big\|_2^2 + \mu\|D\bm{z}_n\|_2^2,$$
where $0 \le \mu \le 1, D = \bm{J}_{n}\bm{Q}_2$ as was discussed in the previous subsection. However, note that the choice of $D$ and $\mu$ depends on a priori knowledge about power lines in the network. Thus an algorithm designer might want to incorporate information, that $X/R$ ratio for a cable is close to 1, by making $D = [1~ -1]$.

\textbf{X/R modification}. One method to improve algorithm performance with respect to noise is to use reactance-to-resistance ratio ($k_n$ for $n$-th power line) usually available from cable specifications. The matrix $\bm{J}_n$ (first defined in (\ref{cost_function-lin})) with $X/R$ modification becomes a vector of size $M$ by $1$, therefore its condition number $\kappa(\bm{J}_n) = 1$, which leads to a better noise resistivity for impedance estimation:
\begin{equation*}
\bm{J}_n = \Real{(\bm{j}_n)} - k_n\Imag{(\bm{j}_n)}.
\end{equation*}
Whenever the algorithm has information about $X/R$ ratio of a cable we add subscript $(\cdot)_{XR}$ to its name, i.e. BCI$_{XR}$ or LBCI$_{XR}$.

\textbf{Decentralisation}. The LBCI/BCI algorithm estimates power line impedances for all nodes by processing them sequentially. Thus, for large scale networks large memory required in order to store information from all smart meters. However, every iteration uses data only from neighbouring nodes, therefore one iteration can be encapsulated within one smart metering device capable of obtaining data from adjacent nodes, i.e. from neighbouring smart meters. This makes the grid estimation procedure secure and scalable.

\begin{figure}[t]
\centering
\begin{tikzpicture}
\draw[dotted, thick] (0,0) -- (1,0);
\node at (4,0) [rectangle,draw] (SM_n) {$\mathbf{[DBCI_{XR}]}_{n}$};
\draw[thick] (1,0) -- (SM_n);
\draw[arrows = {-Stealth[inset=0pt, angle=30:10pt]}, thick] (1,0) -- (2, 0)  node[below]{$\bm{j}_{n}$} node[above]{$z_{n}$};
\draw[thick] (SM_n) -- (7,0);
\draw[dotted, thick] (7,0) -- (8,0);
\draw[arrows = {-Stealth[inset=0pt, angle=30:10pt]}, thick] (SM_n) -- (6,0)  node[below]{$\bm{j}_{n+1}$} node[above]{$z_{n+1}$};
%\draw[arrows = {-Stealth[fill=none, inset=0pt, angle=90:10pt]}, thick] (SM_n) -- (4,-1);
\draw (1,1) -- (3,1) -- (SM_n);
\draw[arrows = {-Stealth[inset=0pt, angle=30:10pt]}] (3, 1) -- (2, 1) node[above]{$\bm{j}_{n}$};
\draw (SM_n) -- (5,1) -- (7,1);
\draw[arrows = {-Stealth[inset=0pt, angle=30:10pt]}] (7, 1) -- (6, 1) node[above]{$\bm{j}_{n+1}$};
%\draw (SM_n) -- (5,-1) -- (7,-1) node[above]{output} node[pos = 1/2, below]{$\bm{j}_{n}, z_{n}$};
\end{tikzpicture}
\caption{Smart meter at $n$-th node. Concept of decentralised BCI with known $X/R$ ratio}
\label{SM_n_bci}
\end{figure}
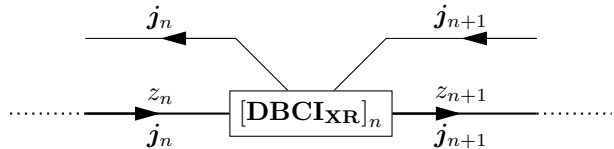

Consider the case of the BCI algorithm and suppose that every smart meter can communicate with its immediate neighbours and the smart meter at node $n+1$ has already calculated values of $\bm{j}_{n+1}$ as in Fig. \ref{SM_n_bci}. Then it sends data to the adjacent smart meter at the $n$-th node, which now has enough information to calculate the impedance of the power line connecting these two meters by performing several iterations of the original BCI algorithm. In addition, it calculates the input current $\mathbf{j}_{n}$ and sends it further up the feeder to the next smart meter and so on. Furthermore, this approach allows each device to have its own information about the $X/R$ ratio for the power lines to which it is connected. The resulting algorithm is called Decentralised BCI with $X/R$ modification, i.e. DBCI$_{XR}$.

\section{Numerical simulations}
\label{num_sim_section}
In order to test the algorithms a single phase power network chain feeder as in Fig. \ref{f1} was implemented in MATLAB/Simulink with $N = 10$. Each power line is modelled as a series connection of resistance and inductance, which is realistic for low voltage overhead lines up to 1-2 km in length \cite{glover2012power}. Loads are modelled as a parallel connection of resistance, inductance and capacitance. Their load shapes are taken from IEEE Low Voltage European Test Feeder data with minute resolution measurements. Although some domestic loads consist of a very complicated composition of resistive, inductive and non-linear parts \cite{barker2013empirical, pipattanasomporn2014load} we limit our simulations as described above, since our algorithm has full information about load impedances via smart meter measurements. We conduct 5000 measurements for all tests, which is equivalent to around 4 days of a network operation.

To illustrate important features of algorithms developed in the paper we consider two separate cases:
\begin{itemize}
  \item No measurement noise. All lines are of the same lengths (50m and 500m).
  \item With measurement noise added (1\%, 0.5\% and 0.1\%). All lines are of the same lengths (50m and 500m).
\end{itemize}
Long lines are chosen in order to show the important properties of the BCI algorithm.

\subsection{LBCI/BCI testing under ideal conditions}
\begin{figure*}[t]
\begin{subfigure}{0.33\linewidth}
\centering
\begin{adjustbox}{width=1\columnwidth}
\begin{tikzpicture}

\begin{axis}[%
width=4.521in,
height=2.8in,
at={(0.758in,0.428in)},
scale only axis,
xmin=1,
xmax=10,
xlabel={Power line number},
xmajorgrids,
ymin=-10,
ymax=0,
ylabel={$\log_{10}\frac{\|z_n - \hat{z}_n\|_2}{\|z_n\|_2}$},
ymajorgrids,
axis background/.style={fill=white},
axis x line*=bottom,
axis y line*=left,
legend style={at={(0.5,1.03)},anchor=south,legend cell align=left,align=left,draw=white!15!black},
xlabel style={font=\LARGE},ylabel style={font=\Huge},legend style={font=\LARGE},ticklabel style={font=\LARGE}
]
\addplot [color=black,dotted,line width=1.0pt,mark=o,mark options={solid}]
  table[row sep=crcr]{%
1	-1.32915348006689\\
2	-1.28838246544068\\
3	-1.26512606312045\\
4	-1.21536180362297\\
5	-1.13471057657888\\
6	-1.1261945142969\\
7	-1.0519885476118\\
8	-1.00224806304511\\
9	-0.92036333225634\\
10	-0.817139363952442\\
};
\addlegendentry{LBCI$_{XR}$};

\addplot [color=red,dotted,line width=1.0pt,mark=triangle,mark options={solid}]
  table[row sep=crcr]{%
1	-3.88090534345849\\
2	-3.86964700560886\\
3	-3.89331250305145\\
4	-3.8735781693953\\
5	-3.80841453476826\\
6	-3.9087658497897\\
7	-3.85832000821522\\
8	-3.84983152727818\\
9	-3.7284474352639\\
10	-3.53956216259043\\
};
\addlegendentry{BCI$_{XR}$ - 1 iteration (LBCI$_{XR}$ - old)};

\addplot [color=blue,dashdotted,line width=1.1pt,mark=triangle,mark options={solid}]
  table[row sep=crcr]{%
1	-4.2904174938234\\
2	-4.2782761821511\\
3	-4.3017398517945\\
4	-4.28214763461237\\
5	-4.21739012037145\\
6	-4.31717221726016\\
7	-4.26705948823946\\
8	-4.25880740157296\\
9	-4.1379507097733\\
10	-3.94946741399456\\
};
\addlegendentry{BCI$_{XR}$ - 10 iterations};

\addplot [color=black,dashed,line width=1.2pt,mark=triangle,mark options={solid}]
  table[row sep=crcr]{%
1	-6.32939409637337\\
2	-6.1831584219856\\
3	-6.18189088210376\\
4	-6.17946507361347\\
5	-6.16732048573507\\
6	-6.1952229482674\\
7	-6.18612232114627\\
8	-6.21054085497848\\
9	-6.16700097486536\\
10	-6.03675702928428\\
};
\addlegendentry{BCI$_{XR}$ - 100 iterations};

\end{axis}
\end{tikzpicture}%
\end{adjustbox}
\end{subfigure}
\begin{subfigure}{0.33\linewidth}
\centering
\begin{adjustbox}{width=1\columnwidth}
\begin{tikzpicture}

\begin{axis}[%
width=4.521in,
height=2.78in,
at={(0.758in,0.426in)},
scale only axis,
xmin=1,
xmax=10,
xlabel={Power line number},
xmajorgrids,
ymin=-10,
ymax=0,
ylabel={$\log_{10}\frac{\|z_n - \hat{z}_n\|_2}{\|z_n\|_2}$},
ymajorgrids,
axis background/.style={fill=white},
axis x line*=bottom,
axis y line*=left,
legend style={at={(0.5,1.03)},anchor=south,legend cell align=left,align=left,draw=white!15!black},
xlabel style={font=\LARGE},ylabel style={font=\Huge},legend style={font=\LARGE},ticklabel style={font=\LARGE}
]
\addplot [color=black,dotted,line width=1.0pt,mark=o,mark options={solid}]
  table[row sep=crcr]{%
1	-1.13177450249518\\
2	-1.17857232776077\\
3	-1.12729593343788\\
4	-1.08402522530668\\
5	-1.08559823954618\\
6	-0.96567971180553\\
7	-0.92846518739287\\
8	-0.797322672080543\\
9	-0.748880390548875\\
10	-0.814977925085405\\
};
\addlegendentry{LBCI};

\addplot [color=red,dotted,line width=1.3pt,mark=triangle,mark options={solid}]
  table[row sep=crcr]{%
1	-3.38652769511604\\
2	-3.3125017754448\\
3	-3.42454108436879\\
4	-3.34459142874997\\
5	-3.21636110316632\\
6	-3.45741188077818\\
7	-3.35822040183589\\
8	-3.37239749939876\\
9	-3.14154963894484\\
10	-2.98392482454151\\
};
\addlegendentry{BCI - 1 iteration (LBCI - old)};

\addplot [color=blue,dashdotted,line width=1.3pt,mark=triangle,mark options={solid}]
  table[row sep=crcr]{%
1	-3.79460707876873\\
2	-3.72099482480146\\
3	-3.83273674153348\\
4	-3.7531238584518\\
5	-3.62532359196459\\
6	-3.86566185875066\\
7	-3.76694141604543\\
8	-3.78111340125635\\
9	-3.55086554012878\\
10	-3.39352482393243\\
};
\addlegendentry{BCI - 10 iterations};

\addplot [color=black,dashed,line width=1.3pt,mark=triangle,mark options={solid}]
  table[row sep=crcr]{%
1	-5.94810422550278\\
2	-5.9401901905404\\
3	-5.92476973786823\\
4	-5.90395647396985\\
5	-5.87782690940461\\
6	-5.8722849057396\\
7	-5.84927456959633\\
8	-5.8329346642505\\
9	-5.7542133357923\\
10	-5.72594239175156\\
};
\addlegendentry{BCI - 100 iterations};

\end{axis}
\end{tikzpicture}%
\end{adjustbox}
\end{subfigure}
\begin{subfigure}{0.33\linewidth}
\centering
\begin{adjustbox}{width=1\columnwidth}
\begin{tikzpicture}

\begin{axis}[%
width=4.521in,
height=3.219in,
at={(0.758in,0.434in)},
scale only axis,
xmin=1,
xmax=10,
xlabel={Power line number},
xmajorgrids,
ymin=0,
ymax=50,
ylabel={$\kappa$(\mbox{\boldmath$J$}$_n$)},
ymajorgrids,
axis background/.style={fill=white},
legend style={at={(0.5,1.03)},anchor=south,legend cell align=left,align=left,draw=white!15!black},
xlabel style={font=\LARGE},ylabel style={font=\Huge},legend style={font=\LARGE},ticklabel style={font=\LARGE}
]
\addplot [color=black,line width=1.0pt,only marks,mark=o,mark options={solid}]
  table[row sep=crcr]{%
1	4.64523681642281\\
2	4.34031112589818\\
3	4.261847277491\\
4	4.02560140250221\\
5	3.58180441458692\\
6	3.69702267447058\\
7	3.34817777291077\\
8	3.36985419389401\\
9	3.02993190395761\\
10	2.35959548143542\\
};
\addlegendentry{Condition number};

\end{axis}
\end{tikzpicture}%
\end{adjustbox}
\end{subfigure}
\caption{Impedance identification error for each power line. FBCI and BCI algorithms under ideal conditions (left, centre). Condition number for matrix $\bm{J}_n$(right).}
\label{z_rel_err}
\end{figure*}
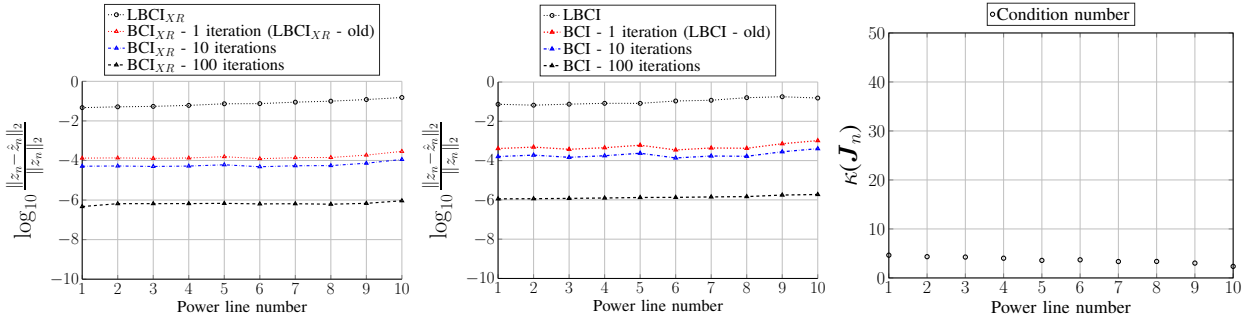
We first show BCI performance under ideal conditions, when measuring devices do not have measurement error. To demonstrate the main features of a proposed algorithm we use oversimplified conditions. Each line has the following parameters: length - 50m, $R = 0.4 ~\Omega/km$, $X/R$ ratio - $0.7$. According to IEEE Low Voltage European Test Feeder data description the power factor is kept constant for each load. This is a good assumption for a high level analysis, however, for our algorithms the more detailed behaviour is required. To model uncertainty in power factor we vary it between 0.9 and 1 as a random Gaussian variable ($\mu = 0.95, \sigma = 0.05$) for each load.

We first show algorithm performance with known $X/R$ ratio, since in this case, the result does not depend on the condition number of the measurement matrix ($\bm{J}_n$). Figure \ref{z_rel_err} (left) illustrates identification relative error of BCI algorithm after different number of iterations. The red dotted line on the figure corresponds to the linearised case, i.e. when we assume a small phase difference ($\Delta_n \ll 1$) between two nodes. Note, that it coincides with a very common in the literature LBCI approach. Blue, and black lines show the relative error after 10 and 100 iterations correspondingly. The step size $\alpha$ was chosen to be 0.1. After about 30 iterations the error is close to simulation accuracy. Note that the same result is unachievable with LBCI or LBCI-old algorithms.

The next figure (Fig. \ref{z_rel_err} - centre) compares BCI algorithm performance under ideal conditions with LBCI. Note, that the corresponding line current matrix $\bm{J}_n$ has a condition number around 4-5 for a given set of load shapes from LV European Test Feeder data (Fig. \ref{z_rel_err} - right) and with power factor disturbances introduced. Thus the performance of algorithms is slightly worse than in the previous case, however the result of BCI algorithm is still very close to the simulation accuracy limit. Again, LBCI-old coincides with the first iteration of the BCI algorithm.

\begin{figure*}[t]
\begin{subfigure}{0.33\linewidth}
\centering
\begin{adjustbox}{width=1\columnwidth}
\begin{tikzpicture}

\begin{axis}[%
width=4.521in,
height=2.493in,
at={(0.758in,0.434in)},
scale only axis,
xmin=0,
xmax=5000,
xtick={   0, 1000, 2000, 3000, 4000, 5000},
xlabel={Measurement number},
xmajorgrids,
ymode=log,
ymin=1,
ymax=1000,
yminorticks=true,
ylabel={$\frac{\|\mathbf{z} - \hat{\mathbf{z}}\|_2}{\|\mathbf{z}\|_2}$, \%},
ymajorgrids,
yminorgrids,
axis background/.style={fill=white},
axis x line*=bottom,
axis y line*=left,
legend style={at={(0.5,1.03)},anchor=south,legend cell align=left,align=left,draw=white!15!black},
xlabel style={font=\LARGE},ylabel style={font=\LARGE},legend style={font=\LARGE},ticklabel style={font=\LARGE}
]
\addplot [color=red,dotted,line width=1.0pt,mark size=0.7pt,mark=*,mark options={solid}]
  table[row sep=crcr]{%
100	618.565874809085\\
200	275.107927604509\\
300	202.986394427141\\
400	174.143346467306\\
500	77.93772915443\\
600	72.4416225257139\\
700	68.130107378111\\
800	65.0667628821093\\
900	63.7561722528446\\
1000	51.7593180051441\\
1100	47.3092775501452\\
1200	42.7948503279262\\
1300	41.5649703657043\\
1400	40.206573276376\\
1500	40.1134762264486\\
1600	40.0608865335942\\
1700	40.1136928732485\\
1800	40.173151501695\\
1900	39.1517458543583\\
2000	37.7875296434527\\
2100	36.5566350133897\\
2200	36.4270020595124\\
2300	36.3712161584117\\
2400	36.2440405611231\\
2500	35.9310860897663\\
2600	35.7248430457915\\
2700	35.5989950715627\\
2800	35.4277719465521\\
2900	35.4012206514555\\
3000	35.3804908540011\\
3100	35.3975097979703\\
3200	35.4492057891983\\
3300	35.4581223995724\\
3400	34.5083499897861\\
3500	34.4032826196327\\
3600	34.3648505927901\\
3700	34.3034650004817\\
3800	34.3232839133271\\
3900	34.2287106974776\\
4000	34.181243284833\\
4100	34.0879353992075\\
4200	33.9460640567723\\
4300	33.823089793748\\
4400	33.8199531414247\\
4500	33.7977770801437\\
4600	33.7683935769202\\
4700	33.7532492472821\\
4800	33.8010418043093\\
4900	33.6587794219313\\
5000	33.6306326494405\\
};
\addlegendentry{LBCI - 1\% full scale error};

\addplot [color=blue,dotted,line width=1.0pt,mark size=0.7pt,mark=*,mark options={solid}]
  table[row sep=crcr]{%
100	352.382634256233\\
200	141.486595016329\\
300	109.528793228904\\
400	93.8828649025291\\
500	47.5474668051383\\
600	43.6459362623189\\
700	41.7148773664865\\
800	40.2174833440735\\
900	39.7397686350978\\
1000	37.6573198651872\\
1100	35.3880567887947\\
1200	33.5605349248709\\
1300	32.7347667333502\\
1400	32.3719595766183\\
1500	32.3609412618906\\
1600	32.3367446608151\\
1700	32.3032423781999\\
1800	32.2632910792546\\
1900	32.1606565516506\\
2000	31.6662257856626\\
2100	31.4059364576298\\
2200	31.3181308351929\\
2300	31.2621034199713\\
2400	31.2897254020552\\
2500	31.2111344695258\\
2600	31.2222036000879\\
2700	31.1827690000903\\
2800	31.1497071445986\\
2900	31.1359988180254\\
3000	31.1321275844053\\
3100	31.141145544122\\
3200	31.1473150797466\\
3300	31.1185537155983\\
3400	30.9393761140801\\
3500	30.8373209837428\\
3600	30.8211447532472\\
3700	30.7965717335581\\
3800	30.7844949112658\\
3900	30.7505007424251\\
4000	30.7411402576624\\
4100	30.7399469762595\\
4200	30.7265848120366\\
4300	30.6692474690435\\
4400	30.6695555195663\\
4500	30.6658082373681\\
4600	30.6653088449481\\
4700	30.6646150175997\\
4800	30.6580565496481\\
4900	30.5782424698311\\
5000	30.5738902999195\\
};
\addlegendentry{LBCI - 0.5\% full scale error};

\addplot [color=black,dotted,line width=1.0pt,mark size=0.7pt,mark=*,mark options={solid}]
  table[row sep=crcr]{%
100	69.6455909615255\\
200	41.3148846012258\\
300	36.8130380065506\\
400	35.2075785676292\\
500	29.8549616548759\\
600	29.244258437695\\
700	29.6032842895404\\
800	29.4773332571922\\
900	29.4976766733221\\
1000	29.4920714755628\\
1100	29.3214686254113\\
1200	29.2570577496857\\
1300	29.1827965922184\\
1400	29.0815350040724\\
1500	29.1063535464245\\
1600	29.1084431680269\\
1700	29.1097548731771\\
1800	29.1073153388602\\
1900	29.1909177452964\\
2000	29.2874136347419\\
2100	29.3544777363983\\
2200	29.3607030663155\\
2300	29.3391808035276\\
2400	29.3621378841458\\
2500	29.3505453721207\\
2600	29.3792440764553\\
2700	29.3841195666729\\
2800	29.3752047242743\\
2900	29.3772157040799\\
3000	29.3786861938124\\
3100	29.3784519248367\\
3200	29.3800834837762\\
3300	29.3744658217687\\
3400	29.5027487784969\\
3500	29.5085088000744\\
3600	29.4990539170773\\
3700	29.5035727521978\\
3800	29.4999150721819\\
3900	29.4943933743683\\
4000	29.5041986213115\\
4100	29.4956716273405\\
4200	29.4984608897722\\
4300	29.460671820908\\
4400	29.4614602370312\\
4500	29.460461105279\\
4600	29.4595629204515\\
4700	29.4589453970515\\
4800	29.4563173683792\\
4900	29.4183972736971\\
5000	29.418409849176\\
};
\addlegendentry{LBCI - 0.1\% full scale error};

\addplot [color=red,solid,line width=1.0pt]
  table[row sep=crcr]{%
100	10939.4282794858\\
200	4218.79293382921\\
300	3281.33474638765\\
400	2599.76617799588\\
500	1170.89270892909\\
600	946.600522161743\\
700	832.457314486713\\
800	784.261008542446\\
900	764.901751111227\\
1000	604.08558306595\\
1100	558.460378328972\\
1200	446.674763923221\\
1300	420.315666155738\\
1400	408.609916813595\\
1500	406.24694564777\\
1600	406.574718286192\\
1700	404.942586687761\\
1800	404.423220903262\\
1900	384.465064708311\\
2000	368.070237505923\\
2100	343.030592031864\\
2200	340.612338475662\\
2300	335.327360270985\\
2400	331.652956219918\\
2500	326.826969927445\\
2600	322.933240020542\\
2700	316.194945360981\\
2800	313.676427429031\\
2900	310.011180861667\\
3000	310.141336859814\\
3100	308.570295534127\\
3200	308.716064436462\\
3300	307.42875518688\\
3400	274.8955869407\\
3500	264.132671945233\\
3600	262.119603904118\\
3700	261.396169546429\\
3800	259.571925345053\\
3900	258.32997816703\\
4000	255.150373500367\\
4100	250.501676227442\\
4200	245.907508667846\\
4300	238.767654757424\\
4400	238.504241236432\\
4500	238.65073691125\\
4600	238.449172471361\\
4700	239.149365510132\\
4800	238.593003038806\\
4900	232.695534502624\\
5000	232.857524441911\\
};
\addlegendentry{BCI - 1\% full scale error (LBCI-old)};

\addplot [color=blue,solid,line width=1.0pt]
  table[row sep=crcr]{%
100	4988.13502781159\\
200	2444.87548208943\\
300	1733.3249857583\\
400	1362.27101754373\\
500	618.336586986953\\
600	531.146403811036\\
700	442.1948107434\\
800	422.191221463147\\
900	403.771786737491\\
1000	314.391176663283\\
1100	286.839637612527\\
1200	239.871235609107\\
1300	225.589302036281\\
1400	212.501253742105\\
1500	211.119466862647\\
1600	210.266081627266\\
1700	209.283377716353\\
1800	208.10170705769\\
1900	191.574028473684\\
2000	180.074443135741\\
2100	168.811777599424\\
2200	167.450518993585\\
2300	165.847538490874\\
2400	165.78422627458\\
2500	163.866224667952\\
2600	160.363277531053\\
2700	157.615577690254\\
2800	155.845124205837\\
2900	154.575970904311\\
3000	154.076091355224\\
3100	153.344013534164\\
3200	153.422976528203\\
3300	153.323579773637\\
3400	139.473256237771\\
3500	135.071348446532\\
3600	134.157554469797\\
3700	133.805378471932\\
3800	133.503207284249\\
3900	133.233338480496\\
4000	130.617395736445\\
4100	130.55064061357\\
4200	129.629425163305\\
4300	127.704815594137\\
4400	127.436836417746\\
4500	127.340412103804\\
4600	127.296991658512\\
4700	126.844783224533\\
4800	126.183289403131\\
4900	122.112369308204\\
5000	121.948639664765\\
};
\addlegendentry{BCI - 0.5\% full scale error (LBCI-old)};

\addplot [color=black,solid,line width=1.0pt]
  table[row sep=crcr]{%
100	983.097984365482\\
200	459.231063318821\\
300	376.804443565766\\
400	302.00130790956\\
500	118.555779963691\\
600	103.3726771843\\
700	86.5674726763889\\
800	82.8863800904153\\
900	81.7801433011419\\
1000	63.9209324169423\\
1100	59.3539145125219\\
1200	44.9947925703019\\
1300	42.8654739846427\\
1400	42.3791929465445\\
1500	42.2467192621202\\
1600	41.9679996512012\\
1700	41.9195047726181\\
1800	41.6675180370154\\
1900	38.4540225521658\\
2000	36.7329981605019\\
2100	34.7268959936169\\
2200	34.2448428725205\\
2300	33.8342804788993\\
2400	33.7059753690994\\
2500	32.9789583380315\\
2600	32.4221252657376\\
2700	32.3104756296501\\
2800	31.9133046882613\\
2900	31.7772443881123\\
3000	31.7701208348717\\
3100	31.7616112489574\\
3200	31.6857654558978\\
3300	31.5960338000104\\
3400	27.8164634410411\\
3500	27.0915026283469\\
3600	26.8150551735293\\
3700	26.7427908700357\\
3800	26.7904096262388\\
3900	26.5361393411003\\
4000	25.9229387582616\\
4100	25.7398585711523\\
4200	25.3460023141886\\
4300	24.9996313706155\\
4400	24.9843284663807\\
4500	24.9646864664521\\
4600	24.9482414111492\\
4700	24.9808542291075\\
4800	25.0062385390108\\
4900	24.5412230366541\\
5000	24.4458111154257\\
};
\addlegendentry{BCI - 0.1\% full scale error (LBCI-old)};

\end{axis}
\end{tikzpicture}%
\end{adjustbox}
\end{subfigure}
\begin{subfigure}{0.33\linewidth}
\centering
\begin{adjustbox}{width=1\columnwidth}
\begin{tikzpicture}

\begin{axis}[%
width=4.521in,
height=2.493in,
at={(0.758in,0.434in)},
scale only axis,
xmin=0,
xmax=5000,
xtick={   0, 1000, 2000, 3000, 4000, 5000},
xlabel={Measurement number},
xmajorgrids,
ymode=log,
ymin=1,
ymax=1000,
yminorticks=true,
ylabel={$\frac{\|\mathbf{z} - \hat{\mathbf{z}}\|_2}{\|\mathbf{z}\|_2}$, \%},
ymajorgrids,
yminorgrids,
axis background/.style={fill=white},
axis x line*=bottom,
axis y line*=left,
legend style={at={(0.5,1.03)},anchor=south,legend cell align=left,align=left,draw=white!15!black},
xlabel style={font=\LARGE},ylabel style={font=\LARGE},legend style={font=\LARGE},ticklabel style={font=\LARGE}
]
\addplot [color=red,dotted,line width=1.0pt,mark size=0.7pt,mark=*,mark options={solid}]
  table[row sep=crcr]{%
100	577.391750400882\\
200	267.093230969695\\
300	205.319242503762\\
400	166.928295300879\\
500	66.644882813857\\
600	61.6679189789651\\
700	55.3852532392064\\
800	50.9604169261538\\
900	49.3675170342762\\
1000	42.0212254866031\\
1100	36.3700843526153\\
1200	29.9193523111009\\
1300	27.2278430482018\\
1400	25.7649706685327\\
1500	25.5327647984282\\
1600	25.6474224322936\\
1700	25.6574706552111\\
1800	25.4680309409088\\
1900	24.5091774128577\\
2000	23.076946555035\\
2100	21.4235929636326\\
2200	21.2421839472294\\
2300	21.1764258356123\\
2400	20.9478765989289\\
2500	20.3576609733448\\
2600	20.2064694206182\\
2700	19.9846228279037\\
2800	19.7587681214098\\
2900	19.6852826688951\\
3000	19.6657931119346\\
3100	19.6560772678979\\
3200	19.739210201449\\
3300	19.6835829537958\\
3400	18.6071014044608\\
3500	17.7729543095237\\
3600	17.7295345972129\\
3700	17.7539063172693\\
3800	17.6991400368571\\
3900	17.6828713873297\\
4000	17.4070259473985\\
4100	17.4172552304756\\
4200	17.1499873281785\\
4300	16.8575389581397\\
4400	16.868921028951\\
4500	16.8657273675102\\
4600	16.8760096189548\\
4700	16.8776479981479\\
4800	16.8635141639981\\
4900	16.497290253911\\
5000	16.5219252269178\\
};
\addlegendentry{LBCI$_{XR}$ - 1\% full scale error};

\addplot [color=blue,dotted,line width=1.0pt,mark size=0.7pt,mark=*,mark options={solid}]
  table[row sep=crcr]{%
100	336.645848518894\\
200	128.554202881726\\
300	107.91341470573\\
400	89.7008586225277\\
500	35.8305282561023\\
600	31.5147132341474\\
700	27.8205188450343\\
800	26.0399639208426\\
900	25.373910928077\\
1000	21.5908272677598\\
1100	19.665633555213\\
1200	16.7821611898977\\
1300	15.7642012447506\\
1400	15.4424522531692\\
1500	15.3629412812489\\
1600	15.3227855127073\\
1700	15.3042926181891\\
1800	15.2809154453315\\
1900	14.9187553765836\\
2000	14.3766246539883\\
2100	13.7706769744652\\
2200	13.6016477410138\\
2300	13.5306006602\\
2400	13.4899507193898\\
2500	13.3924228474751\\
2600	13.3308271756306\\
2700	13.2275946373813\\
2800	13.1932006033215\\
2900	13.2074808783891\\
3000	13.2153492814913\\
3100	13.2121605727348\\
3200	13.2095297553954\\
3300	13.1620105508672\\
3400	12.3088178147556\\
3500	12.1544657262373\\
3600	12.0160245126537\\
3700	11.9270114174846\\
3800	11.8946433629559\\
3900	11.8375568689779\\
4000	11.7292897638418\\
4100	11.6841083720192\\
4200	11.709783417035\\
4300	11.6190952873433\\
4400	11.6143501006531\\
4500	11.592644405011\\
4600	11.5742624064494\\
4700	11.554518002781\\
4800	11.5694556420955\\
4900	11.423074472196\\
5000	11.4175127490441\\
};
\addlegendentry{LBCI$_{XR}$ - 0.5\% full scale error};

\addplot [color=black,dotted,line width=1.0pt,mark size=0.7pt,mark=*,mark options={solid}]
  table[row sep=crcr]{%
100	60.6232061565378\\
200	27.2289448057352\\
300	21.7268329589992\\
400	19.1635426659609\\
500	10.9431204082108\\
600	10.1947024844337\\
700	10.0140096348267\\
800	9.85809285195655\\
900	9.80178458581397\\
1000	9.46487266176593\\
1100	9.28310922112982\\
1200	9.13196395568197\\
1300	9.06044621194058\\
1400	8.93904581145372\\
1500	8.94388026207321\\
1600	8.94197512933743\\
1700	8.94582080975783\\
1800	8.94689282566839\\
1900	8.96661443409027\\
2000	8.99318668403021\\
2100	8.97789317908973\\
2200	8.98827643088499\\
2300	8.9771309821604\\
2400	8.98788038408225\\
2500	8.99255464487493\\
2600	8.98933336210284\\
2700	8.98233202696932\\
2800	8.97437072482077\\
2900	8.98022976246306\\
3000	8.98257483406135\\
3100	8.97941213292436\\
3200	8.97433050709443\\
3300	8.96974470624438\\
3400	9.02449448128123\\
3500	9.02438520888396\\
3600	9.01713476981925\\
3700	9.0206951925146\\
3800	9.01921761024206\\
3900	9.01760644990102\\
4000	9.01844924824081\\
4100	9.00894158839334\\
4200	9.00713181772295\\
4300	8.97885211895225\\
4400	8.97795544281009\\
4500	8.97601619209\\
4600	8.97468330161962\\
4700	8.97311973683783\\
4800	8.97175386011031\\
4900	8.94679929163184\\
5000	8.94983106818811\\
};
\addlegendentry{LBCI$_{XR}$ - 0.1\% full scale error};

\addplot [color=red,solid,line width=1.0pt]
  table[row sep=crcr]{%
100	636.500144257714\\
200	292.676529266655\\
300	223.795513847441\\
400	182.346181281471\\
500	72.2391461623452\\
600	66.5647500460818\\
700	59.8727928669973\\
800	54.9590912238281\\
900	53.0709276860564\\
1000	45.2094835288144\\
1100	38.7612198601652\\
1200	31.6039078982571\\
1300	28.6101670171299\\
1400	27.0686490062978\\
1500	26.8127530888437\\
1600	26.9539956248889\\
1700	26.952166385095\\
1800	26.7900874503672\\
1900	25.8048499970827\\
2000	24.0780007369917\\
2100	22.1563269154779\\
2200	21.9291605690581\\
2300	21.82142918688\\
2400	21.5762728588356\\
2500	20.9036726252275\\
2600	20.7442032531047\\
2700	20.450458762919\\
2800	20.1757940868614\\
2900	20.0890369152345\\
3000	20.0641691491575\\
3100	20.0629870030643\\
3200	20.1631784593289\\
3300	20.1094427897021\\
3400	19.0207477509391\\
3500	18.0006920405027\\
3600	17.9482466160831\\
3700	17.9855167035162\\
3800	17.9127744576634\\
3900	17.8743480688399\\
4000	17.5283624638278\\
4100	17.543397666076\\
4200	17.2208643332019\\
4300	16.8681730969016\\
4400	16.8778226724616\\
4500	16.87658311737\\
4600	16.8880211110168\\
4700	16.8939169017106\\
4800	16.8741869033947\\
4900	16.4524265610638\\
5000	16.4710227186601\\
};
\addlegendentry{BCI$_{XR}$ - 1\% full scale error (LBCI$_{XR}$ - old)};

\addplot [color=blue,solid,line width=1.0pt]
  table[row sep=crcr]{%
100	371.34684258334\\
200	141.004301187621\\
300	117.709262864169\\
400	97.9572088257102\\
500	37.9426917051369\\
600	33.0808380413626\\
700	28.9504901330387\\
800	26.8795429042067\\
900	26.0117379109377\\
1000	21.7736645458502\\
1100	19.3172803174777\\
1200	15.6951061027285\\
1300	14.364146130252\\
1400	13.9415796632522\\
1500	13.8479745269059\\
1600	13.791722702049\\
1700	13.7857650049981\\
1800	13.7599606682894\\
1900	13.2075572092958\\
2000	12.5249694456453\\
2100	11.6182954572368\\
2200	11.395892596789\\
2300	11.3223075802414\\
2400	11.2836252705669\\
2500	11.1598111517771\\
2600	11.0436340217559\\
2700	10.9057020441624\\
2800	10.8514999045513\\
2900	10.8446097573783\\
3000	10.8596901588358\\
3100	10.8484737656465\\
3200	10.861178873008\\
3300	10.7914180988687\\
3400	9.54348935647062\\
3500	9.27100120291306\\
3600	9.07527673743194\\
3700	8.91996663317877\\
3800	8.88739559824935\\
3900	8.78998380310235\\
4000	8.63453034733242\\
4100	8.57720408146672\\
4200	8.61412132044649\\
4300	8.48513726562341\\
4400	8.47561462705886\\
4500	8.44670705081155\\
4600	8.41733050783347\\
4700	8.39226701405558\\
4800	8.42067572309058\\
4900	8.15614971251391\\
5000	8.13779609290959\\
};
\addlegendentry{BCI$_{XR}$ - 0.5\% full scale error (LBCI$_{XR}$ - old)};

\addplot [color=black,solid,line width=1.0pt]
  table[row sep=crcr]{%
100	66.0141848328862\\
200	28.1320284190169\\
300	21.5494579753043\\
400	18.4927117515729\\
500	7.3930446912482\\
600	6.43355363669567\\
700	5.54636554178825\\
800	5.27268540100493\\
900	5.01692690104795\\
1000	4.13232294071725\\
1100	3.80109360438199\\
1200	3.17910427037647\\
1300	3.01203361217865\\
1400	2.77510632773033\\
1500	2.75677356937227\\
1600	2.7531736562209\\
1700	2.75223068924941\\
1800	2.75386088602477\\
1900	2.64407600640884\\
2000	2.46433665573116\\
2100	2.30158312768054\\
2200	2.27918385424576\\
2300	2.25211115169549\\
2400	2.22548639306507\\
2500	2.221428322728\\
2600	2.1676407010861\\
2700	2.13066887267363\\
2800	2.10949762752772\\
2900	2.10159633004459\\
3000	2.10565987103901\\
3100	2.0959100868024\\
3200	2.09131624004766\\
3300	2.08504415791899\\
3400	1.97707125929211\\
3500	1.91400306597603\\
3600	1.89868012631728\\
3700	1.88293101620837\\
3800	1.87764129602836\\
3900	1.88245095796468\\
4000	1.84775182684586\\
4100	1.817556606915\\
4200	1.78064496623178\\
4300	1.73535891599895\\
4400	1.73344057879927\\
4500	1.73023002720965\\
4600	1.72951855580784\\
4700	1.72615288625156\\
4800	1.72459825576531\\
4900	1.70165967800584\\
5000	1.69530490101971\\
};
\addlegendentry{BCI$_{XR}$ - 0.1\% full scale error (LBCI$_{XR}$ - old)};

\end{axis}
\end{tikzpicture}%
\end{adjustbox}
\end{subfigure}
\begin{subfigure}{0.33\linewidth}
\centering
\begin{adjustbox}{width=1\columnwidth}
\input{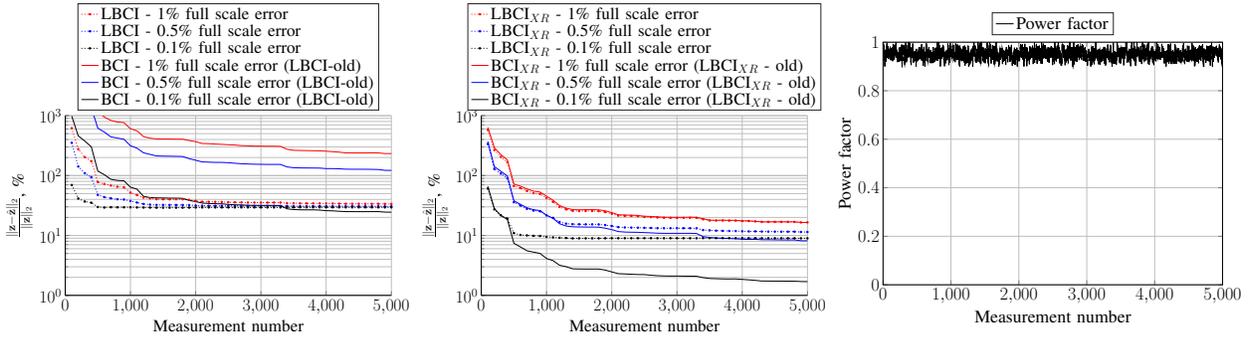}
\end{adjustbox}
\end{subfigure}
\caption{Dependence of impedance identification error on measurement number. BCI algorithm under noisy conditions. Averaged over 100 realisations.}
\label{z_rel_err3}
\end{figure*}

\begin{figure*}[t]
\begin{subfigure}{0.33\linewidth}
\centering
\begin{adjustbox}{width=1\columnwidth}
\begin{tikzpicture}

\begin{axis}[%
width=4.521in,
height=2.493in,
at={(0.758in,0.434in)},
scale only axis,
xmin=0,
xmax=5000,
xtick={   0, 1000, 2000, 3000, 4000, 5000},
xlabel={Measurement number},
xmajorgrids,
ymode=log,
ymin=1,
ymax=1000,
yminorticks=true,
ylabel={$\frac{\|\mathbf{z} - \hat{\mathbf{z}}\|_2}{\|\mathbf{z}\|_2}$, \%},
ymajorgrids,
yminorgrids,
axis background/.style={fill=white},
axis x line*=bottom,
axis y line*=left,
legend style={at={(0.5,1.03)},anchor=south,legend cell align=left,align=left,draw=white!15!black},
xlabel style={font=\LARGE},ylabel style={font=\LARGE},legend style={font=\LARGE},ticklabel style={font=\LARGE}
]
\addplot [color=red,dotted,line width=1.0pt,mark size=0.7pt,mark=*,mark options={solid}]
  table[row sep=crcr]{%
100	638.949193481054\\
200	285.70632182257\\
300	209.502339854108\\
400	171.638806921245\\
500	56.645261210907\\
600	54.22432967775\\
700	48.4070108550094\\
800	45.5579631046186\\
900	43.9265013147077\\
1000	39.6900200621635\\
1100	37.3139967910915\\
1200	33.7202706964587\\
1300	32.2708677517544\\
1400	30.8519306407674\\
1500	30.6299745994681\\
1600	30.4741274960733\\
1700	30.4411214252344\\
1800	30.4765370150899\\
1900	29.0628210361043\\
2000	26.4273624315949\\
2100	24.4719832958088\\
2200	23.7458857941618\\
2300	23.5066633297105\\
2400	22.9717915741404\\
2500	22.5405706381273\\
2600	22.1226461839415\\
2700	22.0209644511223\\
2800	21.8561817997409\\
2900	21.8300263708918\\
3000	21.7992148647942\\
3100	21.8198750420406\\
3200	21.8975662766937\\
3300	21.9280484672161\\
3400	20.5566598250784\\
3500	20.5709226407569\\
3600	20.5750150741861\\
3700	20.5289081253015\\
3800	20.5841350806246\\
3900	20.3059484817058\\
4000	20.0412008661088\\
4100	19.7646851911023\\
4200	19.7248735737716\\
4300	19.6916455272616\\
4400	19.6390876030231\\
4500	19.6102477112483\\
4600	19.5229409113847\\
4700	19.4589219229\\
4800	19.4507801909066\\
4900	18.8390541766195\\
5000	18.6845331429238\\
};
\addlegendentry{LBCI - 1\% full scale error};

\addplot [color=blue,dotted,line width=1.0pt,mark size=0.7pt,mark=*,mark options={solid}]
  table[row sep=crcr]{%
100	365.948825420284\\
200	148.216446076282\\
300	112.82942755769\\
400	92.6660881577041\\
500	40.1555990892812\\
600	38.8016657747387\\
700	34.362425175556\\
800	31.9028270853374\\
900	30.9455346352107\\
1000	28.0279998990404\\
1100	26.053396529023\\
1200	24.1610712800741\\
1300	23.9164814349801\\
1400	23.2776564430657\\
1500	23.1007184966155\\
1600	22.9496166699504\\
1700	22.8789805796042\\
1800	22.8170153791079\\
1900	21.8639913285255\\
2000	19.477866340094\\
2100	17.9227592133653\\
2200	17.1375749638645\\
2300	16.8926394652853\\
2400	16.5999628494894\\
2500	16.3933460694156\\
2600	16.1966786578198\\
2700	16.0862140966498\\
2800	16.0565275103474\\
2900	16.0020972844761\\
3000	15.9774325053264\\
3100	15.9694056087035\\
3200	15.9742574063977\\
3300	15.9453833135101\\
3400	15.4657982888896\\
3500	15.3954911357374\\
3600	15.4315000387376\\
3700	15.4603364027782\\
3800	15.4576722731547\\
3900	15.3018064174461\\
4000	15.0390875342421\\
4100	14.8501237659851\\
4200	15.0801491176738\\
4300	15.1959547847747\\
4400	15.1534624825223\\
4500	15.1401647686803\\
4600	15.0874400903603\\
4700	15.029286036236\\
4800	14.9565018376587\\
4900	14.325497747777\\
5000	14.1980293385826\\
};
\addlegendentry{LBCI - 0.5\% full scale error};

\addplot [color=black,dotted,line width=1.0pt,mark size=0.7pt,mark=*,mark options={solid}]
  table[row sep=crcr]{%
100	74.4598834190606\\
200	44.4950743040255\\
300	34.8342380330379\\
400	28.7987073250332\\
500	32.1936194582918\\
600	31.9308237946469\\
700	27.9299784433982\\
800	26.0440102073967\\
900	25.328942043142\\
1000	22.077479335062\\
1100	20.6430801399838\\
1200	19.6388477630566\\
1300	19.7070865082712\\
1400	19.139488227619\\
1500	18.9882857765568\\
1600	18.8551996439041\\
1700	18.8153744934514\\
1800	18.7979378498159\\
1900	17.8512708972745\\
2000	15.6959684149122\\
2100	14.3415927643204\\
2200	13.6530117410485\\
2300	13.4285744595613\\
2400	13.224660233667\\
2500	13.1144657761498\\
2600	12.9189365193015\\
2700	12.8116020666747\\
2800	12.8059658557099\\
2900	12.741825074648\\
3000	12.730315460627\\
3100	12.7014160326678\\
3200	12.6958870656686\\
3300	12.7000573127662\\
3400	12.6847533920869\\
3500	12.7920922286813\\
3600	12.8356589584009\\
3700	12.9141790563003\\
3800	12.9371800409294\\
3900	12.786776324107\\
4000	12.4876544206742\\
4100	12.2599982874187\\
4200	12.5216743294418\\
4300	12.7023157184054\\
4400	12.6564170330362\\
4500	12.6506955677774\\
4600	12.6000811421575\\
4700	12.5400349844273\\
4800	12.4662420549318\\
4900	11.8096026473679\\
5000	11.6825148899311\\
};
\addlegendentry{LBCI - 0.1\% full scale error};

\addplot [color=red,solid,line width=1.0pt]
  table[row sep=crcr]{%
100	22728.0844401917\\
200	5118.69657000351\\
300	3536.99850914306\\
400	816.769790166386\\
500	217.789962578717\\
600	186.673358272156\\
700	170.588712436484\\
800	158.672266693795\\
900	154.945990194417\\
1000	102.812960987144\\
1100	90.0903675033844\\
1200	79.8235730543928\\
1300	76.0952811801202\\
1400	73.2939001233863\\
1500	73.3105008058256\\
1600	72.7745763727444\\
1700	72.7783815896182\\
1800	72.718031543639\\
1900	70.6605592162223\\
2000	66.0800688049927\\
2100	63.0546347055478\\
2200	62.2194878973759\\
2300	62.0956124039623\\
2400	60.6097719919569\\
2500	60.0059486729373\\
2600	59.3725361804331\\
2700	59.273596157382\\
2800	58.9970650417395\\
2900	59.2004854022767\\
3000	59.2540643803915\\
3100	59.3824432036232\\
3200	59.5270961745916\\
3300	59.5457679559992\\
3400	57.0456483739896\\
3500	56.5590202126003\\
3600	56.3618548265959\\
3700	56.3310515498867\\
3800	56.3326145456274\\
3900	56.0118499285963\\
4000	55.4987848173866\\
4100	55.3352463267593\\
4200	54.7210962583346\\
4300	53.8589408369959\\
4400	53.6297037181513\\
4500	53.6839536356162\\
4600	53.2751694592901\\
4700	53.2522253135632\\
4800	53.177736128386\\
4900	52.5670815106881\\
5000	52.3268073688193\\
};
\addlegendentry{BCI - 1\% full scale error (LBCI-old)};

\addplot [color=blue,solid,line width=1.0pt]
  table[row sep=crcr]{%
100	11400.6317104915\\
200	2709.11003762981\\
300	1837.480807797\\
400	460.735052239631\\
500	102.005394167206\\
600	87.2906130185545\\
700	80.0915493937165\\
800	76.8472619532403\\
900	74.1173835937968\\
1000	51.9491088048463\\
1100	44.8245611687652\\
1200	40.5784875779692\\
1300	37.5024040866603\\
1400	35.7216733221953\\
1500	35.7193452101075\\
1600	35.5593681721333\\
1700	35.5168223438536\\
1800	35.5311590526479\\
1900	34.7140095542001\\
2000	31.6135326665586\\
2100	30.1820496229043\\
2200	29.251150624789\\
2300	28.9584491275208\\
2400	27.7120685560186\\
2500	27.2715114539494\\
2600	27.0170947009927\\
2700	26.9242291975569\\
2800	26.8678826767531\\
2900	26.8440375438859\\
3000	26.772186808124\\
3100	26.7609659888137\\
3200	26.7669739770353\\
3300	26.7189769781261\\
3400	26.6314209655037\\
3500	26.4145706764439\\
3600	26.281180847016\\
3700	26.18937750428\\
3800	26.1581314055083\\
3900	26.086201997588\\
4000	26.0601240884831\\
4100	25.9459986670809\\
4200	25.6442362748963\\
4300	25.5423909233786\\
4400	25.6431195654694\\
4500	25.5931984861965\\
4600	25.5147322817117\\
4700	25.4599795184458\\
4800	25.4167488602764\\
4900	25.4113724076065\\
5000	25.3137359321453\\
};
\addlegendentry{BCI - 0.5\% full scale error (LBCI-old)};

\addplot [color=black,solid,line width=1.0pt]
  table[row sep=crcr]{%
100	2426.97935227751\\
200	520.097429495387\\
300	370.627433664178\\
400	86.4088147076159\\
500	21.4660738897338\\
600	17.6791178477347\\
700	15.9634979839217\\
800	14.8705928411826\\
900	14.3726130409558\\
1000	10.8271106697378\\
1100	9.40872402129834\\
1200	8.06864024640509\\
1300	7.34031789104499\\
1400	7.06138614146183\\
1500	7.05484216150573\\
1600	7.0294981712566\\
1700	7.01212772537375\\
1800	7.00416054756797\\
1900	6.82498779889169\\
2000	6.48510124758661\\
2100	6.15065684701161\\
2200	6.00827285737952\\
2300	5.97754929288951\\
2400	5.9144943049684\\
2500	5.89186642858549\\
2600	5.85979970044402\\
2700	5.83257088974041\\
2800	5.82381641579081\\
2900	5.8059225153119\\
3000	5.79328015943364\\
3100	5.78378652643428\\
3200	5.778141003319\\
3300	5.76138406166361\\
3400	5.56721982716132\\
3500	5.51920557437916\\
3600	5.50204025657614\\
3700	5.47477625359662\\
3800	5.46857327003434\\
3900	5.43009535981615\\
4000	5.40738490087627\\
4100	5.37831097108823\\
4200	5.28639416421409\\
4300	5.27038746815866\\
4400	5.2594982138701\\
4500	5.25765191351608\\
4600	5.23462274053984\\
4700	5.2121104596335\\
4800	5.22313758218916\\
4900	5.1755789476994\\
5000	5.1637319100943\\
};
\addlegendentry{BCI - 0.1\% full scale error (LBCI-old)};

\end{axis}
\end{tikzpicture}%
\end{adjustbox}
\end{subfigure}
\begin{subfigure}{0.33\linewidth}
\centering
\begin{adjustbox}{width=1\columnwidth}
\begin{tikzpicture}

\begin{axis}[%
width=4.521in,
height=2.493in,
at={(0.758in,0.434in)},
scale only axis,
xmin=0,
xmax=5000,
xtick={   0, 1000, 2000, 3000, 4000, 5000},
xlabel={Measurement number},
xmajorgrids,
ymode=log,
ymin=1,
ymax=1000,
yminorticks=true,
ylabel={$\frac{\|\mathbf{z} - \hat{\mathbf{z}}\|_2}{\|\mathbf{z}\|_2}$, \%},
ymajorgrids,
yminorgrids,
axis background/.style={fill=white},
axis x line*=bottom,
axis y line*=left,
legend style={at={(0.5,1.03)},anchor=south,legend cell align=left,align=left,draw=white!15!black},
xlabel style={font=\LARGE},ylabel style={font=\LARGE},legend style={font=\LARGE},ticklabel style={font=\LARGE}
]
\addplot [color=red,dotted,line width=1.0pt,mark size=0.7pt,mark=*,mark options={solid}]
  table[row sep=crcr]{%
100	564.481923572343\\
200	272.484646683434\\
300	208.426318945265\\
400	158.213670454915\\
500	43.2530423217083\\
600	41.5182412254465\\
700	38.9661632298446\\
800	35.8667829694198\\
900	34.6453567193563\\
1000	33.1857591055094\\
1100	31.1448505693434\\
1200	27.9294712318906\\
1300	26.3580797869923\\
1400	25.2236408437764\\
1500	24.9756551590762\\
1600	25.0349164236565\\
1700	25.0712026123535\\
1800	24.886564853032\\
1900	23.9468278266095\\
2000	22.7014139042234\\
2100	21.2480197058099\\
2200	21.0848721371781\\
2300	20.9468040563461\\
2400	20.4269987914412\\
2500	19.8025716616455\\
2600	19.5042199594777\\
2700	19.2122754282362\\
2800	19.0071690887514\\
2900	18.935948398865\\
3000	18.9203080737734\\
3100	18.9171750201994\\
3200	18.9611118308726\\
3300	18.9235284080245\\
3400	18.042555271127\\
3500	17.3935526311081\\
3600	17.3417739884359\\
3700	17.3415275499434\\
3800	17.2870452840383\\
3900	17.2456765872908\\
4000	16.8401098422829\\
4100	16.8408153178072\\
4200	16.6176224272498\\
4300	16.3100484947216\\
4400	16.3282373326122\\
4500	16.3218421831683\\
4600	16.3411522434454\\
4700	16.3294078129004\\
4800	16.3008558277689\\
4900	15.7063308788233\\
5000	15.7273364097281\\
};
\addlegendentry{LBCI$_{XR}$ - 1\% full scale error};

\addplot [color=blue,dotted,line width=1.0pt,mark size=0.7pt,mark=*,mark options={solid}]
  table[row sep=crcr]{%
100	329.925863089506\\
200	131.225400166045\\
300	108.682116612618\\
400	83.7199061694903\\
500	28.8230690447383\\
600	27.7650730482484\\
700	25.3653666029142\\
800	24.1420469821693\\
900	23.319025328123\\
1000	21.3910778175977\\
1100	20.239584858566\\
1200	18.5791291104949\\
1300	18.0789314845483\\
1400	17.5735396804153\\
1500	17.3909669079807\\
1600	17.3511084937124\\
1700	17.2916814906363\\
1800	17.2592605367503\\
1900	16.6891412816758\\
2000	15.8888483539762\\
2100	15.0836693636405\\
2200	14.8880404487268\\
2300	14.6889304580055\\
2400	14.6656772990501\\
2500	14.3434279926162\\
2600	14.1022012744251\\
2700	13.8509685795753\\
2800	13.748250227589\\
2900	13.7058469895653\\
3000	13.7142507230477\\
3100	13.695916710193\\
3200	13.6834915722055\\
3300	13.6562095902141\\
3400	12.7693963651823\\
3500	12.5799418412041\\
3600	12.4092738524808\\
3700	12.3274627923577\\
3800	12.2783319083166\\
3900	12.1352816496976\\
4000	11.8427028028654\\
4100	11.7446357855642\\
4200	11.7806664465127\\
4300	11.6677949714233\\
4400	11.6707793429469\\
4500	11.6511645275251\\
4600	11.6402943665509\\
4700	11.6130716498999\\
4800	11.5847218495694\\
4900	11.2462644521625\\
5000	11.2210900998497\\
};
\addlegendentry{LBCI$_{XR}$ - 0.5\% full scale error};

\addplot [color=black,dotted,line width=1.0pt,mark size=0.7pt,mark=*,mark options={solid}]
  table[row sep=crcr]{%
100	60.7260001960621\\
200	28.5136917518507\\
300	21.776027499514\\
400	17.9259430391606\\
500	20.9443736945084\\
600	20.8758614150486\\
700	19.3364298631374\\
800	18.6803834934567\\
900	18.147328216363\\
1000	16.3891053836993\\
1100	15.4959649300077\\
1200	14.674399916599\\
1300	14.6404061425242\\
1400	14.0457178688981\\
1500	13.9048430723567\\
1600	13.856960474772\\
1700	13.8093330493652\\
1800	13.7762342416552\\
1900	13.0387989871401\\
2000	12.2797144693579\\
2100	11.5930283912766\\
2200	11.5616790113835\\
2300	11.403174069441\\
2400	11.5720258998527\\
2500	11.3178370570721\\
2600	11.0260163491567\\
2700	10.7779347486385\\
2800	10.6592282753327\\
2900	10.6051121440769\\
3000	10.6051026734001\\
3100	10.5756435746973\\
3200	10.5599332659085\\
3300	10.5504080690842\\
3400	10.0816336064368\\
3500	9.94396168508698\\
3600	9.86586674722158\\
3700	9.84590160006396\\
3800	9.83298772494378\\
3900	9.7237746939626\\
4000	9.46287268266297\\
4100	9.39366871198762\\
4200	9.40094153031319\\
4300	9.32162687665783\\
4400	9.33169172397399\\
4500	9.32441675808094\\
4600	9.3414630748198\\
4700	9.33513432309673\\
4800	9.28491438338282\\
4900	8.95135307849097\\
5000	8.91797640279913\\
};
\addlegendentry{LBCI$_{XR}$ - 0.1\% full scale error};

\addplot [color=red,solid,line width=1.0pt]
  table[row sep=crcr]{%
100	687.162632312096\\
200	310.581909851338\\
300	238.8564667258\\
400	180.902080862412\\
500	53.4955991216033\\
600	51.2362052694262\\
700	47.625529793001\\
800	43.4151969883932\\
900	41.6042688251828\\
1000	39.1592652954947\\
1100	35.6674632497747\\
1200	30.8666349120025\\
1300	28.4308514040445\\
1400	27.1421605790363\\
1500	26.8893304853057\\
1600	27.0088309173177\\
1700	27.0210911435296\\
1800	26.9021906597949\\
1900	25.909629026642\\
2000	24.2939280348576\\
2100	22.2330831591075\\
2200	21.9430333743929\\
2300	21.7602911615194\\
2400	21.0996757855984\\
2500	20.4152038530171\\
2600	20.1781764715397\\
2700	19.8531190698737\\
2800	19.6471796084147\\
2900	19.5783091473384\\
3000	19.5478853043168\\
3100	19.5685040159804\\
3200	19.6250956269494\\
3300	19.5835760920241\\
3400	18.8266934184721\\
3500	17.9666396640691\\
3600	17.9064595928047\\
3700	17.9216792159022\\
3800	17.8410373739834\\
3900	17.7969598729175\\
4000	17.3979067603842\\
4100	17.4239973806169\\
4200	17.1328370985482\\
4300	16.7613593265288\\
4400	16.7802145012919\\
4500	16.7797938421224\\
4600	16.7930000447739\\
4700	16.7899473869768\\
4800	16.7567820457399\\
4900	16.1031637272547\\
5000	16.1241967612528\\
};
\addlegendentry{BCI$_{XR}$ - 1\% full scale error (LBCI$_{XR}$-old)};

\addplot [color=blue,solid,line width=1.0pt]
  table[row sep=crcr]{%
100	404.1321725797\\
200	149.907145372373\\
300	125.047013550393\\
400	96.2052957147922\\
500	28.1385540829226\\
600	26.2564948383921\\
700	23.5119455154061\\
800	21.9201928831837\\
900	20.9847690368135\\
1000	19.282089588279\\
1100	17.6908447294148\\
1200	15.1549445078433\\
1300	13.8749438847117\\
1400	13.4550865853242\\
1500	13.3511605321021\\
1600	13.2977969667917\\
1700	13.3066388691057\\
1800	13.2644462935719\\
1900	12.8539181640328\\
2000	12.4206026227885\\
2100	11.5917373605801\\
2200	11.3306951547318\\
2300	11.2297673687889\\
2400	11.1164398891805\\
2500	10.9186750128686\\
2600	10.7415289122971\\
2700	10.6389535941284\\
2800	10.6098105589105\\
2900	10.5977162124967\\
3000	10.6166732692131\\
3100	10.5963378544181\\
3200	10.6086783421659\\
3300	10.5527427481681\\
3400	9.51993532860792\\
3500	9.27804611209175\\
3600	9.10848054090651\\
3700	8.96508734171332\\
3800	8.91796616030325\\
3900	8.79864730528896\\
4000	8.64131244856672\\
4100	8.5756305258051\\
4200	8.61643132577958\\
4300	8.49587863563992\\
4400	8.48170845152912\\
4500	8.45965768252504\\
4600	8.42125014536799\\
4700	8.38815775564816\\
4800	8.41531848948757\\
4900	8.10499014595432\\
5000	8.0932767770725\\
};
\addlegendentry{BCI$_{XR}$ - 0.5\% full scale error (LBCI$_{XR}$-old)};

\addplot [color=black,solid,line width=1.0pt]
  table[row sep=crcr]{%
100	72.7195490887134\\
200	29.8355381810924\\
300	22.9326955663768\\
400	18.3217119954536\\
500	5.35976617652142\\
600	4.9693423911885\\
700	4.38776740234361\\
800	4.16555527409481\\
900	3.93519315037906\\
1000	3.62564767137419\\
1100	3.47519256870061\\
1200	3.04084066537432\\
1300	2.93004121813041\\
1400	2.72205052231501\\
1500	2.70965158616435\\
1600	2.70543137930573\\
1700	2.70071097718423\\
1800	2.69858839141467\\
1900	2.59940062923238\\
2000	2.45141671070456\\
2100	2.2937077370773\\
2200	2.24764004009169\\
2300	2.20854708761487\\
2400	2.1515491404239\\
2500	2.14573839185848\\
2600	2.08475565438333\\
2700	2.05198197165185\\
2800	2.03458243419119\\
2900	2.02860411526808\\
3000	2.03171625846985\\
3100	2.02131144427194\\
3200	2.01716151247194\\
3300	2.01121763815109\\
3400	1.9230171551157\\
3500	1.87301500214421\\
3600	1.85730571413849\\
3700	1.84271459294255\\
3800	1.83633521278918\\
3900	1.84035199739431\\
4000	1.79664021587933\\
4100	1.76738193462632\\
4200	1.73275273990007\\
4300	1.68679283100553\\
4400	1.68237342382349\\
4500	1.67865029296265\\
4600	1.67946483079034\\
4700	1.6775715993085\\
4800	1.67466385165876\\
4900	1.64625904730214\\
5000	1.64095509018859\\
};
\addlegendentry{BCI$_{XR}$ - 0.1\% full scale error (LBCI$_{XR}$ - old)};

\end{axis}
\end{tikzpicture}%
\end{adjustbox}
\end{subfigure}
\begin{subfigure}{0.33\linewidth}
\centering
\begin{adjustbox}{width=1\columnwidth}
\input{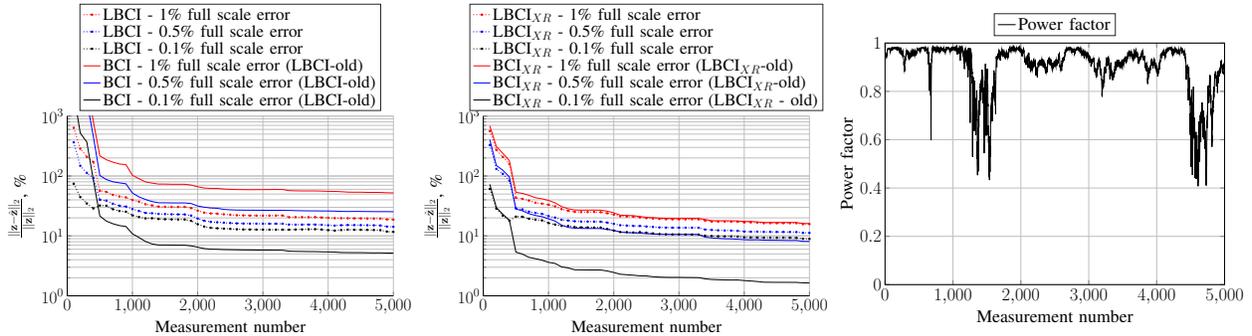}
\end{adjustbox}
\end{subfigure}
\caption{Dependence of impedance identification error on measurement number. BCI algorithm under noisy conditions. Averaged over 100 realisations.}
\label{z_rel_err4}
\end{figure*}

\subsection{LBCI/BCI testing with modifications under noisy conditions}
In this subsection we introduce Gaussian noise to model different accuracy classes (1\%FS, 0.5\%FS and 0.1\%FS) of commercially available smart meters and test FBCI and BCI algorithms performance under these conditions. Usually, accuracy class is specified as a percentage of a full scale with 95\% confidence interval (i.e. two standard deviations in our case). Again, we conduct two tests: when $X/R$ ratio is given and without it. In addition we consider a scenario when power factor variations are much larger, i.e. Fig \ref{z_rel_err4} (right) as an example for such variations for a particular load. 

Figure \ref{z_rel_err3} shows the identification error dependence on the measurements number for three accuracy classes listed above using LBCI and BCI algorithms. In Figure \ref{z_rel_err3} (left) the results are shown when $X/R$ ratio is unknown a priori. In this case LBCI outperform both BCI and conventional LBCI-old algorithms, providing more accurate identification as well as significantly smaller amount of computations. However, note how BCI algorithm outperforms LBCI in the Figure \ref{z_rel_err3} (right), i.e. when condition number of $\kappa(\bm{J}_n) = 1$. This shows that when the condition number of matrix $\bm{J}_n$ decreases, the quality of BCI estimation significantly increases. The same applies to LBCI-old algorithm, its estimation results almost coincide with BCI ones. To demonstrate the effect of condition number of $\bm{J}_n$, we conducted an additional simulation, using the same feeder, but with much higher power factor variations (see Figure \ref{z_rel_err4}). Note, how BCI algorithm outperforms LBCI for measurements with $0.1\%$ accuracy by a factor of 4 when $XR$ ratio is unknown, and by a factor of 5 when it is known, almost achieving $1\%$ of impedances estimation accuracy. 
\begin{remark}
Although, LBCI algorithm shows better results for $1\%$ and $0.5\%$ measurement accuracy cases there are two key properties should be mentioned:
\begin{itemize}
\item The gap between BCI and LBCI performance curves for different power factor variation scenarios is reduced by a factor of 5 and 3 respectively, which shows that BCI algorithms outperforms LBCI provided enough variation in the grid. 
\item BCI performance can achieve the same (even slightly better according to corollary for Proposition \ref{prop_bci_vs_lbci}) level of accuracy by adding a regularisation term.
\end{itemize}
\end{remark}

\subsection{BCI testing with long power lines}

\begin{figure}[H]
\centering
\begin{adjustbox}{width=0.9\columnwidth}
\begin{tikzpicture}

\begin{axis}[%
width=6.028in,
height=3.769in,
at={(1.011in,0.509in)},
scale only axis,
xmin=0,
xmax=5000,
xtick={   0, 1000, 2000, 3000, 4000, 5000},
xlabel={Measurement number},
xmajorgrids,
ymode=log,
ymin=1,
ymax=1000,
yminorticks=true,
ylabel={$\frac{\|\mathbf{z} - \hat{\mathbf{z}}\|_2}{\|\mathbf{z}\|_2}$, \%},
ymajorgrids,
yminorgrids,
axis background/.style={fill=white},
axis x line*=bottom,
axis y line*=left,
legend style={at={(0.5,1.03)},anchor=south,legend cell align=left,align=left,draw=white!15!black},
xlabel style={font=\LARGE},ylabel style={font=\LARGE},legend style={font=\LARGE},ticklabel style={font=\LARGE}
]
\addplot [color=red,dotted,line width=1.0pt,mark size=0.7pt,mark=*,mark options={solid}]
  table[row sep=crcr]{%
100	659.040534857702\\
200	306.80206106422\\
300	202.922606522803\\
400	184.27915420211\\
500	74.4704195738183\\
600	65.302035719461\\
700	57.3290890863622\\
800	53.0782732352832\\
900	52.323212175496\\
1000	41.2297135264796\\
1100	37.06800933764\\
1200	31.1708657492995\\
1300	28.5307090304024\\
1400	27.44627439343\\
1500	27.0745167846451\\
1600	26.9631606437801\\
1700	26.8984670467172\\
1800	26.9084358644316\\
1900	25.9551686685729\\
2000	25.108726644597\\
2100	22.2956337843983\\
2200	22.2633406784104\\
2300	22.2916275857289\\
2400	21.9397445129273\\
2500	21.7571665685396\\
2600	21.3864445891304\\
2700	20.8868182867011\\
2800	20.8551033191808\\
2900	20.9010994549578\\
3000	20.8686650046051\\
3100	20.9124703000949\\
3200	20.8649528940537\\
3300	20.879100147923\\
3400	20.2736837839125\\
3500	20.1870360368623\\
3600	20.1483733715194\\
3700	20.1434208685339\\
3800	20.1541992817378\\
3900	20.120723953138\\
4000	19.9277625715534\\
4100	19.7542004915796\\
4200	19.6090014861816\\
4300	19.2160115690463\\
4400	19.1669184580374\\
4500	19.1719221538792\\
4600	19.1452710662481\\
4700	19.1159603514197\\
4800	19.0748744457345\\
4900	18.8335536134415\\
5000	18.8003659417871\\
};
\addlegendentry{LBCI - 1\% full scale error};

\addplot [color=blue,dotted,line width=1.0pt,mark size=0.7pt,mark=*,mark options={solid}]
  table[row sep=crcr]{%
100	325.939763095285\\
200	155.305033657335\\
300	110.411672317769\\
400	95.4146877978017\\
500	41.249436770124\\
600	36.0799825249198\\
700	31.5515105083569\\
800	29.7407468237229\\
900	28.5139786412904\\
1000	23.3469491449701\\
1100	21.2283873729015\\
1200	17.5892945136175\\
1300	16.3329394225103\\
1400	15.5986178330847\\
1500	15.3392836535708\\
1600	15.3030886303937\\
1700	15.3013568747385\\
1800	15.2723262564557\\
1900	14.8839722243997\\
2000	14.4851641842197\\
2100	14.1518047703867\\
2200	14.1183832560406\\
2300	14.0758549689415\\
2400	13.9949994416347\\
2500	13.9070829485743\\
2600	13.7552356520581\\
2700	13.505073096113\\
2800	13.3534689842058\\
2900	13.2658197596707\\
3000	13.2278015007836\\
3100	13.2300607379895\\
3200	13.2102023486308\\
3300	13.186987940092\\
3400	12.5409473242492\\
3500	12.3511909319057\\
3600	12.26799573637\\
3700	12.1353927885668\\
3800	12.1142811230124\\
3900	12.0523914594524\\
4000	11.9651291018913\\
4100	11.9173615199099\\
4200	11.7982989393014\\
4300	11.8609818939798\\
4400	11.8495505040594\\
4500	11.8453263970241\\
4600	11.8344695876821\\
4700	11.8186064625387\\
4800	11.7672816465301\\
4900	11.7743649326105\\
5000	11.7632499127752\\
};
\addlegendentry{LBCI - 0.5\% full scale error};

\addplot [color=black,dotted,line width=1.0pt,mark size=0.7pt,mark=*,mark options={solid}]
  table[row sep=crcr]{%
100	71.8022323347998\\
200	33.5281508422428\\
300	23.3159265416903\\
400	20.9351102185999\\
500	13.7629070258335\\
600	11.358207818472\\
700	10.8396584874401\\
800	11.1827066979431\\
900	10.6439571015401\\
1000	10.9752048229656\\
1100	10.4062099477565\\
1200	10.5142233534481\\
1300	9.91349750262337\\
1400	9.59147621152229\\
1500	9.47820241714138\\
1600	9.45209137703198\\
1700	9.42662008795498\\
1800	9.41885346762396\\
1900	9.3923495929566\\
2000	9.52883194479701\\
2100	9.84218757470693\\
2200	9.80857438950287\\
2300	9.80742495264447\\
2400	9.72782068631855\\
2500	9.6449617653382\\
2600	9.53036152927272\\
2700	9.3478429288162\\
2800	9.22839299045763\\
2900	9.1612579858623\\
3000	9.15255463183914\\
3100	9.16315026376884\\
3200	9.16010683777694\\
3300	9.13424602543443\\
3400	9.29776501168745\\
3500	9.18245897227697\\
3600	9.07642578779062\\
3700	9.0649821038938\\
3800	9.06338048113857\\
3900	9.05060245662827\\
4000	9.09696308306439\\
4100	9.02984178062694\\
4200	8.91456351897199\\
4300	8.99263264164221\\
4400	8.98700538504772\\
4500	8.98781087357142\\
4600	8.98245974680313\\
4700	8.96250748953804\\
4800	8.93097907012001\\
4900	9.02799943896461\\
5000	9.00963044088909\\
};
\addlegendentry{LBCI - 0.1\% full scale error};

\addplot [color=red,solid,line width=1.0pt]
  table[row sep=crcr]{%
100	659.56881658749\\
200	307.301241447016\\
300	201.32307985644\\
400	183.048146929529\\
500	71.4734203225816\\
600	62.8829862479565\\
700	56.0648964454405\\
800	51.7403645671444\\
900	50.8909332085476\\
1000	39.5480465087212\\
1100	35.3693373314896\\
1200	29.1246375562055\\
1300	26.1757073825054\\
1400	25.2293016732408\\
1500	24.8535576227902\\
1600	24.742855389003\\
1700	24.6866695156716\\
1800	24.6919392958537\\
1900	23.6626260819765\\
2000	22.5470417196135\\
2100	19.456433102685\\
2200	19.4338467039465\\
2300	19.4249202277284\\
2400	19.0647680565012\\
2500	18.8452931322232\\
2600	18.4315363725663\\
2700	17.998981725382\\
2800	18.0342982314189\\
2900	18.1199328673537\\
3000	18.0830753664742\\
3100	18.1088959258068\\
3200	18.0591309631655\\
3300	18.0740087008344\\
3400	17.2248911077859\\
3500	17.1330901203646\\
3600	17.1559751066928\\
3700	17.1881760649542\\
3800	17.1892395968817\\
3900	17.1431805387667\\
4000	16.9495792721079\\
4100	16.8409564873455\\
4200	16.7311388506264\\
4300	16.2190677222075\\
4400	16.1692003673591\\
4500	16.1819692862289\\
4600	16.1554878603417\\
4700	16.1310605934899\\
4800	16.1226976603156\\
4900	15.8201195044776\\
5000	15.7940658119176\\
};
\addlegendentry{BCI - 1\% full scale error};

\addplot [color=blue,solid,line width=1.0pt]
  table[row sep=crcr]{%
100	327.111291797483\\
200	155.631950217755\\
300	110.152025521629\\
400	95.1843497018664\\
500	38.7358972183351\\
600	33.9790764486301\\
700	29.8055176386943\\
800	27.7758985195908\\
900	26.6045603632109\\
1000	20.6658420954446\\
1100	18.5161824556313\\
1200	14.0668412004198\\
1300	12.8532247560692\\
1400	12.1647364321055\\
1500	11.9698602049376\\
1600	11.9378700347329\\
1700	11.9497633267475\\
1800	11.9080679001962\\
1900	11.4411598126037\\
2000	10.6677937860089\\
2100	9.87080049468944\\
2200	9.83105050952715\\
2300	9.75560992404414\\
2400	9.70910135096232\\
2500	9.67797957182692\\
2600	9.54918274226053\\
2700	9.40447390659514\\
2800	9.33203218866336\\
2900	9.29613939549279\\
3000	9.25475556728678\\
3100	9.24050835944374\\
3200	9.20531206054066\\
3300	9.18604070916017\\
3400	8.17666594618074\\
3500	7.92497093708109\\
3600	7.88840505840538\\
3700	7.71424338983087\\
3800	7.69615959142446\\
3900	7.63045170843426\\
4000	7.45434934133202\\
4100	7.45480988173055\\
4200	7.39458664101081\\
4300	7.36457569856779\\
4400	7.35810291141632\\
4500	7.35347974605921\\
4600	7.34662168331193\\
4700	7.35388079006209\\
4800	7.31580144655943\\
4900	7.24861915723705\\
5000	7.24832820182254\\
};
\addlegendentry{BCI - 0.5\% full scale error};

\addplot [color=black,solid,line width=1.0pt]
  table[row sep=crcr]{%
100	71.9541762646144\\
200	33.5714208423819\\
300	23.0963014618009\\
400	20.6938865825754\\
500	7.99384261411036\\
600	6.92782930806457\\
700	6.25376681653603\\
800	5.58500041996704\\
900	5.3338379867869\\
1000	4.29574030341234\\
1100	3.91052795901051\\
1200	3.01098491176513\\
1300	2.66009076324519\\
1400	2.5432381035151\\
1500	2.54630514270835\\
1600	2.53772286418481\\
1700	2.53798012076776\\
1800	2.53935276467103\\
1900	2.50148722650604\\
2000	2.3327724674815\\
2100	2.03064016710436\\
2200	1.99815988831018\\
2300	1.97475963079394\\
2400	1.96831698852112\\
2500	1.9391794510729\\
2600	1.93326065993609\\
2700	1.91211405745849\\
2800	1.9012560051618\\
2900	1.89898745589482\\
3000	1.89969735235255\\
3100	1.88970325620367\\
3200	1.88874581183903\\
3300	1.8820945576694\\
3400	1.70843521121686\\
3500	1.68882991354695\\
3600	1.68490571962326\\
3700	1.66341663541364\\
3800	1.65638371434164\\
3900	1.65086937753352\\
4000	1.60719972069634\\
4100	1.5945916054125\\
4200	1.5933903638057\\
4300	1.55597196508194\\
4400	1.56023326749214\\
4500	1.56199464407532\\
4600	1.55901604396616\\
4700	1.55944716858258\\
4800	1.56132651300968\\
4900	1.55149656017617\\
5000	1.54661416554074\\
};
\addlegendentry{BCI - 0.1\% full scale error};

\end{axis}
\end{tikzpicture}%
\end{adjustbox}
\caption{Comparison of LBCI-old and BCI algorithms in the case of long power lines}
\label{z_rel_long_lines}
\end{figure}
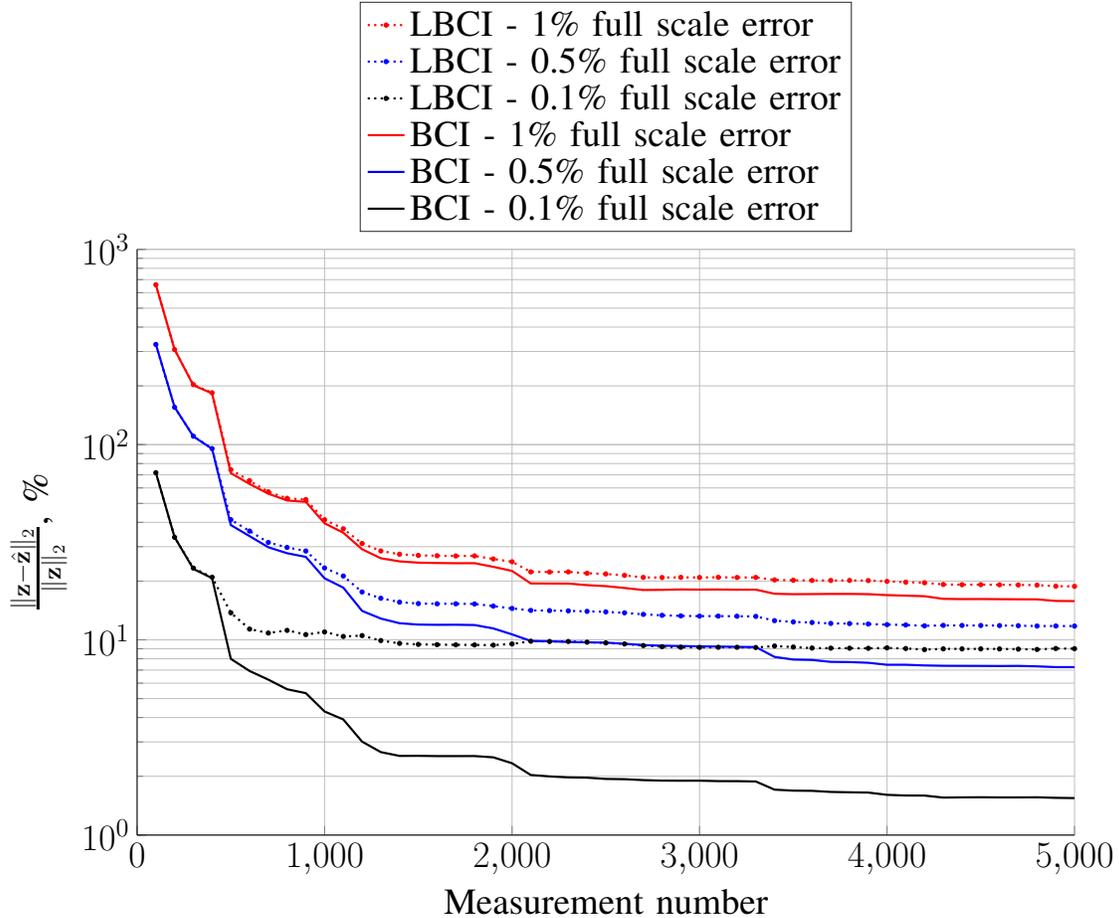

The previous tests under noisy condition did not show the main advantage of using BCI over LBCI-old and, in some cases, LBCI algorithms. However, it can be observed when using long lines in the network, which happens quite often in sparsely populated places. The same simulation was conducted with the lines 500m each with unknown $X/R$ ratios. The results are shown in the Figure \ref{z_rel_long_lines}.

It is noticeable how BCI outperforms LBCI-old by almost 50\% in the case of 1\% of the measurement error, by the factor of 2 when the error is 0.5\% and by the factor of 10 when the measurements are very accurate (0.1\%). Combining this with the previous results let us conclude that BCI algorithm (with regularisation in some cases) is recommended for usage instead of the conventional approach (LBCI-old).

\section{Conclusion}
A novel power line impedance estimation method for the low-voltage grid has been proposed in this paper and a decentralised implementation suitable for smart meters has been developed. The BCI algorithm is iterative and based on least squares estimation making it very attractive for practical applications. We do not require phase synchronous measurements, and we propose modifications that take advantage of known $X/R$ ratio to develop real-time distributed versions of our algorithms. Simulations in MATLAB have shown algorithm performance and key properties for different numbers of measurements and smart meters accuracy classes. 
We have also provided a theoretical framework that allows solution of a certain class of non-convex problems that can be generalised in future work. We will also address possible extensions of the BCI algorithm for topology identification and fault detection.

\section{Acknowledgements}
The authors would like to thank Steven Law and Iman Shames for useful comments and discussions.

\bibliographystyle{IEEEtran}
\bibliography{Refs}

% Generated by IEEEtran.bst, version: 1.14 (2015/08/26)
\begin{thebibliography}{10}
\providecommand{\url}[1]{#1}
\csname url@samestyle\endcsname
\providecommand{\newblock}{\relax}
\providecommand{\bibinfo}[2]{#2}
\providecommand{\BIBentrySTDinterwordspacing}{\spaceskip=0pt\relax}
\providecommand{\BIBentryALTinterwordstretchfactor}{4}
\providecommand{\BIBentryALTinterwordspacing}{\spaceskip=\fontdimen2\font plus
\BIBentryALTinterwordstretchfactor\fontdimen3\font minus
  \fontdimen4\font\relax}
\providecommand{\BIBforeignlanguage}[2]{{%
\expandafter\ifx\csname l@#1\endcsname\relax
\typeout{** WARNING: IEEEtran.bst: No hyphenation pattern has been}%
\typeout{** loaded for the language `#1'. Using the pattern for}%
\typeout{** the default language instead.}%
\else
\language=\csname l@#1\endcsname
\fi
#2}}
\providecommand{\BIBdecl}{\relax}
\BIBdecl

\bibitem{huang2012state}
Y.-F. Huang, S.~Werner, J.~Huang, N.~Kashyap, and V.~Gupta, ``State estimation
  in electric power grids: Meeting new challenges presented by the requirements
  of the future grid,'' \emph{IEEE Signal Processing Magazine}, vol.~29, no.~5,
  pp. 33--43, 2012.

\bibitem{della2014electrical}
D.~Della~Giustina, M.~Pau, P.~A. Pegoraro, F.~Ponci, and S.~Sulis, ``Electrical
  distribution system state estimation: measurement issues and challenges,''
  \emph{IEEE Instrumentation \& Measurement Magazine}, vol.~17, no.~6, pp.
  36--42, 2014.

\bibitem{lam2012distributed}
A.~Y. Lam, B.~Zhang, and N.~T. David, ``Distributed algorithms for optimal
  power flow problem,'' in \emph{2012 IEEE 51st IEEE Conference on Decision and
  Control (CDC)}.\hskip 1em plus 0.5em minus 0.4em\relax IEEE, 2012, pp.
  430--437.

\bibitem{tarkiainen2004identification}
A.~Tarkiainen, R.~Pollanen, M.~Niemela, and J.~Pyrhonen, ``Identification of
  grid impedance for purposes of voltage feedback active filtering,''
  \emph{IEEE Power Electronics Letters}, vol.~2, no.~1, pp. 6--10, 2004.

\bibitem{wood2012power}
A.~J. Wood and B.~F. Wollenberg, \emph{Power generation, operation, and
  control}.\hskip 1em plus 0.5em minus 0.4em\relax John Wiley \& Sons, 2012.

\bibitem{jahangiri2013distributed}
P.~Jahangiri and D.~C. Aliprantis, ``Distributed volt/var control by pv
  inverters,'' \emph{IEEE Transactions on power systems}, vol.~28, no.~3, pp.
  3429--3439, 2013.

\bibitem{de2010synchronized}
J.~De~La~Ree, V.~Centeno, J.~S. Thorp, and A.~G. Phadke, ``Synchronized phasor
  measurement applications in power systems,'' \emph{IEEE Transactions on Smart
  Grid}, vol.~1, no.~1, pp. 20--27, 2010.

\bibitem{yang2010online}
J.~Yang, W.~Li, T.~Chen, W.~Xu, and M.~Wu, ``Online estimation and application
  of power grid impedance matrices based on synchronised phasor measurements,''
  \emph{IET generation, transmission \& distribution}, vol.~4, no.~9, p. 1052,
  2010.

\bibitem{cavraro2015data}
G.~Cavraro, R.~Arghandeh, K.~Poolla, and A.~Von~Meier, ``Data-driven approach
  for distribution network topology detection,'' in \emph{2015 IEEE Power \&
  Energy Society General Meeting}.\hskip 1em plus 0.5em minus 0.4em\relax IEEE,
  2015, pp. 1--5.

\bibitem{deka2016estimating}
D.~Deka, S.~Backhaus, and M.~Chertkov, ``Estimating distribution grid
  topologies: A graphical learning based approach,'' in \emph{Power Systems
  Computation Conference (PSCC), 2016}.\hskip 1em plus 0.5em minus 0.4em\relax
  IEEE, 2016, pp. 1--7.

\bibitem{alahakoon2016smart}
D.~Alahakoon and X.~Yu, ``Smart electricity meter data intelligence for future
  energy systems: A survey,'' \emph{IEEE Transactions on Industrial
  Informatics}, vol.~12, no.~1, pp. 425--436, 2016.

\bibitem{han2016automated}
S.~Han, D.~Kodaira, S.~Han, B.~Kwon, Y.~Hasegawa, and H.~Aki, ``An automated
  impedance estimation method in low-voltage distribution network for
  coordinated voltage regulation,'' \emph{IEEE Transactions on Smart Grid},
  vol.~7, no.~2, pp. 1012--1020, 2016.

\bibitem{cobreces2009grid}
S.~Cobreces, E.~J. Bueno, D.~Pizarro, F.~J. Rodriguez, and F.~Huerta, ``Grid
  impedance monitoring system for distributed power generation electronic
  interfaces,'' \emph{IEEE Transactions on Instrumentation and Measurement},
  vol.~58, no.~9, pp. 3112--3121, 2009.

\bibitem{ciobotaru2011line}
M.~Ciobotaru, V.~Agelidis, and R.~Teodorescu, ``Line impedance estimation using
  model based identification technique,'' in \emph{Power Electronics and
  Applications (EPE 2011), Proceedings of the 2011-14th European Conference
  on}.\hskip 1em plus 0.5em minus 0.4em\relax IEEE, 2011, pp. 1--9.

\bibitem{tariq2016electricity}
M.~Tariq and H.~V. Poor, ``Electricity theft detection and localization in
  grid-tied microgrids,'' \emph{IEEE Transactions on Smart Grid}, 2016.

\bibitem{van2011electricity}
F.~Van Der~Bergh, P.~Kadurek, S.~Cobben, and W.~Kling, ``Electricity theft
  localization based on smart metering,'' in \emph{21st International
  Conference on Electricity Distribution, Frankfurt}, 2011, pp. 1--4.

\bibitem{sahoo2015electricity}
S.~Sahoo, D.~Nikovski, T.~Muso, and K.~Tsuru, ``Electricity theft detection
  using smart meter data,'' in \emph{Innovative Smart Grid Technologies
  Conference (ISGT), 2015 IEEE Power \& Energy Society}.\hskip 1em plus 0.5em
  minus 0.4em\relax IEEE, 2015, pp. 1--5.

\bibitem{diestel2005graph}
R.~Diestel, ``Graph theory, ser,'' \emph{Graduate Texts in Mathematics.
  Springer-Verlag, Heidelberg}, vol. 173, 2005.

\bibitem{SMmin_au}
\BIBentryALTinterwordspacing
B.~R.~W. Group, ``Smart metering infrastructure minimum functionality
  specification,'' ~, Tech. Rep., 2011. [Online]. Available:
  \url{https://link.aemo.com.au/sites/wcl/smartmetering/Pages/BRWG.aspx}
\BIBentrySTDinterwordspacing

\bibitem{NMI_au}
NMI, ``M 6-1 electricity meters, part 1: Metrological and technical
  requirements,'' National Measurement Institute, Australia, Tech. Rep., 2012.

\bibitem{brice1982comparison}
C.~Brice, ``Comparison of approximate and exact voltage drop calculations for
  distribution lines,'' \emph{IEEE Transactions on Power Apparatus and
  Systems}, no.~11, pp. 4428--4431, 1982.

\bibitem{golub1999tikhonov}
G.~H. Golub, P.~C. Hansen, and D.~P. O'Leary, ``Tikhonov regularization and
  total least squares,'' \emph{SIAM Journal on Matrix Analysis and
  Applications}, vol.~21, no.~1, pp. 185--194, 1999.

\bibitem{peppanen2015distribution}
J.~Peppanen, M.~J. Reno, R.~J. Broderick, and S.~Grijalva, ``Distribution
  system secondary circuit parameter estimation for model calibration.'' Sandia
  National Lab.(SNL-NM), Albuquerque, NM (United States), Tech. Rep., 2015.

\bibitem{glover2012power}
J.~D. Glover, M.~S. Sarma, and T.~Overbye, \emph{Power System Analysis \&
  Design, SI Version, chapter 5}.\hskip 1em plus 0.5em minus 0.4em\relax
  Cengage Learning, 2012.

\bibitem{barker2013empirical}
S.~Barker, S.~Kalra, D.~Irwin, and P.~Shenoy, ``Empirical characterization and
  modeling of electrical loads in smart homes,'' in \emph{Green Computing
  Conference (IGCC), 2013 International}.\hskip 1em plus 0.5em minus
  0.4em\relax IEEE, 2013, pp. 1--10.

\bibitem{pipattanasomporn2014load}
M.~Pipattanasomporn, M.~Kuzlu, S.~Rahman, and Y.~Teklu, ``Load profiles of
  selected major household appliances and their demand response
  opportunities,'' \emph{IEEE Transactions on Smart Grid}, vol.~5, no.~2, pp.
  742--750, 2014.

\end{thebibliography}

\section{Appendices}

\appendix
\label{appendix}

\subsection*{Theoretical framework}

Consider the following optimisation problem:
\begin{equation}
\begin{aligned}
& \min_{\bm{x}, \bm{y}} 
& & f(\bm{x}, \bm{y}) \\
& \st
& & \bm{y} = \bm{g}(\bm{x}).
\end{aligned}
\label{min_f_xy}
\end{equation}
where:
$f(\bm{x}, \bm{y}) = \|\bm{A}\bm{x} - \bm{B}\bm{y} - \bm{c} \|_2^2$ with $\bm{A} \in \mathbb{R}^{m\times n}$, $\bm{B} \in \mathbb{R}^{m\times n}$, $\bm{c} \in \mathbb{R}^{m}$ and $\bm{g}(\bm{x})$ is a continuous and differentiable function. Introduce $\bm{h}(\bm{y}) \coloneqq \argmin_{\bm{x}} f(\bm{x}, \bm{y})$

\begin{lemma}
\label{lemma1}
Assume:
\begin{enumerate}[(1)]
\item $\bm{A}$ and $\bm{B}$ are full rank matrices;
\item $\bm{g} \colon \mathbb{R}^m \rightarrow \mathbb{R}^m$ is surjective (onto) mapping; 
\end{enumerate}
Then $\bm{y}^* = \bm{g} \circ \bm{h}(\bm{y}^*)$ and $\bm{x}^* = \bm{h}(\bm{y}^*)$ give the solution $(\bm{x}^*, \bm{y}^*)$ of the optimisation problem (\ref{min_f_xy}).
\end{lemma}

\begin{proof}
$\bm{h}(\bm{y}) = \bm{A}^{\dagger} (\bm{B}\bm{y} + \bm{c})$ is a linear least squares solution. Provided that $\bm{A}$ is a full rank matrix, $\bm{h}(\bm{y})$ exists and unique for any $\bm{y}$. Thus $\bm{h}(\bm{y})$ is an injective mapping.

Next, note that $f(\bm{x}, \bm{y})$ is strictly convex. KKT conditions for (\ref{min_f_xy}) take the form:
\begin{equation}
\begin{split}
{}& \nabla_{\bm{x}} f(\bm{x}, \bm{y}) + \nabla_{\bm{y}} f(\bm{x}, \bm{y}) \nabla_{\bm{x}} \bm{g}(\bm{x}) = 0, \\
& \bm{y} = \bm{g}(\bm{x}).
\end{split}
\end{equation}
Since $\bm{B}$ is a full rank matrix, $\nabla_{\bm{y}} f(\bm{x}, \bm{y}) \Big|_{\bm{x} = \bm{h}(\bm{y})} = \bm{0}$ and therefore $\bm{h}(\bm{y})$ is a unique minimiser of $f(\bm{x}, \bm{y})$. 

Observe that a solution of the problem (\ref{min_f_xy}) will necessarily satisfy:
$
\begin{cases}
\bm{x} = \bm{h}(\bm{y}), \\
\bm{y} = \bm{g}(\bm{x}), 
\end{cases}
$
where from we conclude that the unique $\bm{y}$ can be obtained if $\bm{g} \circ \bm{h}$ form the inverse mapping, i.e. $\bm{h}$ is the right inverse for $\bm{g}$ and $\bm{g}$ is the left inverse for $\bm{h}$. Therefore, using supplementary proposition that are given below:
\begin{itemize}
\item The function is injective if and only if it has a left inverse,
\item The function is surjective if and only if it has a right inverse,
\end{itemize}
we prove the statement of the lemma.
\end{proof}

\subsection*{Supplementary propositions}

\begin{proposition}
Let $X$ and $Y$ are non-empty sets, then $f\colon X \rightarrow Y$ is injective if and only if there exists $g\colon Y \rightarrow X$ such that $g \circ f = id_X$.
\end{proposition}
\begin{proof}
1) Suppose there is $g\colon Y \rightarrow X$ such that $g \circ f = id_X$. Let $f(x) = f(y)$ then $x = id_X(x) = g \circ f(x) = g \circ f(y) = id_X(y) = y$,
i.e. $f(x) = f(y) \Rightarrow x = y$ which means that $f: X \rightarrow Y$ is injective.

2) Suppose $f\colon X \rightarrow Y$ is injective, i.e. $f(x) = f(y) \Rightarrow x = y$. Then if $y \in \image{f}$, there is unique $x \in X: ~f(x) = y$. Fix some $x \in X$ and define $g\colon Y \rightarrow X$ as follows:
$$g(y) = 
\begin{cases} 
f^{-1}(y), & \mbox{if } y \in \image{f}, \\
x, & \mbox{if } y \notin \image{f}, 
\end{cases}$$
Note that $g \circ f (x) = f^{-1}(f(x)) = x$ by injectivity of $f$ and construction of $g$.
\end{proof}

\begin{proposition}
Let $X$ and $Y$ are non-empty sets, then $f\colon X \rightarrow Y$ is surjective if and only if there exists $g\colon Y \rightarrow X$ such that $f \circ g = id_Y$.
\end{proposition}
\begin{proof}
1) Suppose there is $g\colon Y \rightarrow X$ such that $f \circ g = id_Y$. Then for each $y \in Y$ there is $x_y = g(y) \in X$ such that $f \circ g(y) = f(x_y) = id_Y = y$, i.e. f is surjective.

2) Suppose $f\colon X \rightarrow Y$ is surjective. Then for all $y \in Y$ there is $x_y \in X$ such that $f(x_y) = y$. Define $g \colon Y \rightarrow X$ so that it maps each $y$ to $x_y$. Then $\forall y \in Y:~ f \circ g(y) = f(x_y) = y $, 
i.e. $f \circ g = id_Y$.
\end{proof}

\subsection{Proof of Proposition \ref{prop_bci_vs_lbci}}
\begin{proposition*}
Let $(\bm{z}_n^{LBCI}, \bm{1})$ is the solution found by the LBCI algorithm and $(\bm{z}_n^{BCI}, \bm{\gamma}^{BCI}_{n})$ is the solution found by the BCI algorithm.
Also, let $N(\bm{z}_n, \bm{\gamma}_{n}) \coloneqq \big\| \bm{V}_{n-1}\bm{\gamma}_{n} - \bm{v}_{n} - \bm{J}_{n}\bm{Q}_1\bm{z}_{n} \big\|_2^2 + \big\|\bm{V}_{n-1}\sqrt{\bm{1} - \bm{\gamma}^2_{n}} - \bm{J}_{n}\bm{Q}_2\bm{z}_{n} \big\|_2^2$ then
$$N(\bm{z}_n^{BCI}, \bm{\gamma}^{BCI}_{n}) \le N(\bm{z}_n^{LBCI}, \bm{1}),$$
i.e. the BCI algorithm finds the solution that is closer to the optimum of (\ref{pf_network}) than that of LBCI.
\end{proposition*}
\begin{proof}
Consider a linearised BCI algorithm for the power line between nodes $n-1$ and $n$. The solution it finds satisfies:
\begin{equation*}
\begin{split}
{}& N(\bm{z}_n^{LBCI}, \bm{1}) = \big\| \bm{V}_{n-1} - \bm{v}_{n} - \bm{J}_{n}\bm{Q}_1\bm{z}^{LBCI}_{n} \big\|_2^2 + \\
& + \big\|\bm{J}_{n}\bm{Q}_2\bm{z}^{LBCI}_{n} \big\|_2^2 \ge \big\| \bm{V}_{n-1} - \bm{v}_{n} - \bm{J}_{n}\bm{Q}_1\bm{z}^{LBCI}_{n} \big\|_2^2.
\end{split}
\end{equation*}

The corresponding expressions for $\bm{z}_n^{BCI}$ and $\bm{z}_n^{LBCI}$ are given by:  
\begin{equation*}
\begin{split}
{}& \bm{z}_n^{BCI} = \big[ \bm{J}_{n}\bm{Q}_1 \big]^{\dagger} (\bm{V}_{n-1}\bm{\gamma}^{BCI}_{n} - \bm{v}_{n}), \\
& \bm{z}_n^{LBCI} = \big[ \bm{J}_{n}\bm{Q}_1 \big]^{\dagger} (\bm{V}_{n-1} - \bm{v}_{n}),
\end{split}
\end{equation*}
where $\bm{\gamma}^{BCI}_{n} \le \bm{1}$. Therefore:
\begin{equation*}
\begin{split}
{}& N(\bm{z}_n^{BCI}, \bm{\gamma}^{BCI}_{n}) = \big\| \bm{V}_{n-1}\bm{\gamma}^{BCI}_{n} - \bm{v}_{n} - \bm{J}_{n}\bm{Q}_1\bm{z}^{BCI}_{n} \big\|_2^2 \le\\
& \Big\|\big[ \bm{I} - \bm{J}_{n}\bm{Q}_1\big[ \bm{J}_{n}\bm{Q}_1 \big]^{\dagger}\big](\bm{V}_{n-1}\bm{\gamma}^{BCI}_{n} - \bm{v}_{n})\Big\|_2^2 \le \\
& \Big\|\big[ \bm{I} - \bm{J}_{n}\bm{Q}_1\big[ \bm{J}_{n}\bm{Q}_1 \big]^{\dagger}\big](\bm{V}_{n-1} - \bm{v}_{n})\Big\|_2^2 = \\
& \big\| \bm{V}_{n-1} - \bm{v}_{n} - \bm{J}_{n}\bm{Q}_1\bm{z}^{LBCI}_{n} \big\|_2^2 \le N(\bm{z}_n^{LBCI}, \bm{1}).
\end{split}
\end{equation*}
\end{proof}
%\subsection{Proof of Corollary for Proposition \ref{prop_bci_vs_lbci}}
\begin{corollary} In the presence of high measurement noise, we can introduce a regularisation term (as in the LBCI case) for BCI algorithm so that:
$N(\bm{z}_n^{BCI}, \bm{\gamma}^{BCI}_{n}) \le N(\bm{z}_n^{LBCI}, \bm{1})$.
\end{corollary}
\begin{proof}
To see this, note that:
\begin{equation*}
\begin{split}
{}& N(\bm{z}_n^{BCI}, \bm{\gamma}^{BCI}_{n}) \le \Big\|\big[ \bm{I} - \bm{J}_{n}\bm{Q}_1\big[ \bm{J}_{n}\bm{Q}_1 \big]^{\dagger}\big](\bm{V}_{n-1} - \bm{v}_{n})\Big\|_2^2 \le \\
& \Big\|\big[ \bm{I} - \bm{J}_{n}\bm{Q}_1\big[ \bm{J}_{n}\bm{Q}_1 \big]^{\dagger}\big](\bm{V}_{n-1} - \bm{v}_{n})\Big\|_2^2 + \mu \big\|\bm{J}_{n}\bm{Q}_2\bm{z}^{BCI}_{n} \big\|_2^2 \le \\
& \Big\|\big[ \bm{I} - \bm{J}_{n}\bm{Q}_1\big[ \bm{J}_{n}\bm{Q}_1 \big]^{\dagger}\big](\bm{V}_{n-1} - \bm{v}_{n})\Big\|_2^2 + \big\|\bm{J}_{n}\bm{Q}_2\bm{z}^{LBCI}_{n} \big\|_2^2 = \\
& N(\bm{z}_n^{LBCI}, \bm{1}),
\end{split}
\end{equation*}
where $0 \le \mu \le 1$ is a regularisation parameter.
\end{proof}

\end{document}